\documentclass[12pt,twoside]{script_s}
\usepackage{epsfig}
\usepackage{amssymb}
\def\IMSmarkvadjust{0 pt}
\def\IMSmarkhadjust{0 pt}
\def\IMSmarkhpadding{0 pt}
\def\SBIMSMark#1#2#3{
 \font\SBF=cmss10 at 10 true pt
 \font\SBI=cmssi10 at 10 true pt
 \setbox0=\hbox{\SBF \hbox to \IMSmarkhpadding{\relax}
                Stony Brook IMS Preprint \##1}
 \setbox2=\hbox to \wd0{\hfil \SBI #2}
 \setbox4=\hbox to \wd0{\hfil \SBI #3}
 \setbox6=\hbox to \wd0{\hss
             \vbox{\hsize=\wd0 \parskip=0pt \baselineskip=10 true pt
                   \copy0 \break%
                   \copy2 \break%
                   \copy4 \break}}
 \dimen0=\ht6   \advance\dimen0 by \vsize \advance\dimen0 by 8 true pt
                \advance\dimen0 by -\pagetotal
	        \advance\dimen0 by \IMSmarkvadjust
 \dimen2=\hsize \advance\dimen2 by .25 true in
	        \advance\dimen2 by \IMSmarkhadjust

  \openin2=publishd.tex
  \ifeof2\setbox0=\hbox to 0pt{}
  \else 
     \setbox0=\hbox to 3.1 true in{
                \vbox to \ht6{\hsize=3 true in \parskip=0pt  \noindent  
                {\SBI Published in modified form:}\hfil\break
                \input publishd.tex 
                \vfill}}
  \fi
  \closein2
  \ht0=0pt \dp0=0pt
 \ht6=0pt \dp6=0pt
 \setbox8=\vbox to \dimen0{\vfill \hbox to \dimen2{\copy0 \hss \copy6}}
 \ht8=0pt \dp8=0pt \wd8=0pt
 \copy8
 \message{*** Stony Brook IMS Preprint #1, #2. #3 ***}
}

\newtheorem{theo}{Theorem} 
\newtheorem{cor}{Corollary} 
\newtheorem{pro}{Proposition}
\newtheorem{lem}{Lemma}
\newtheorem{defi}{Definition}

\begin{document}
\vspace{-10mm}
\title{Correspondence and Translation Principles for the
Mandelbrot set}
\author{{\large Karsten Keller}\vspace{-1mm}\\
{\large Institut f\"ur Mathematik und Informatik}\vspace{-2mm}\\
{\large Ernst-Moritz-Arndt-Universit\"at}\vspace{-2mm}\\ 
{\large 17487 Greifswald,
Germany}}
\date{}
\thispagestyle{empty}
\maketitle
\renewcommand{\labelenumi}{(\alph{enumi})}\vspace{-15mm}

\begin{abstract}
\noindent New insights into the combinatorial
structure of the
 the Mandelbrot set are given by 
`Correspondence' and 
`Translation' Principles both conjectured and partially proved by E.\ Lau
and D.\ Schleicher. 
We provide complete proofs of these principles and
discuss results related to them.
\end{abstract}
\def\IMSmarkvadjust{-.3in}
\SBIMSMark{1997/14}{October 1997}{Revised April 1999}

\vskip -.2in
\noindent{\sectfont Note:} The Translation and Correspondence Principles stated in this
note do not hold in the general case.  Please see the errata at the end for
details about which parts of the two statements are incorrect and which
parts remain true.

\section{The conjectures by Lau and Schleicher}
\paragraph{Introduction.}
The detailed structure of the Mandelbrot set $M$ 
is extremely comp\-li\-ca\-ted.
However, much of the structure can be described by different kinds
of symmetry and self-similarity. (For a listing of various symmetries 
in $M$, see \cite{schlei6}.)

For example, each neighborhood of a boundary point of $M$ contains
infinitely many topological copies of $M$ itself.
This is a consequence of the (unpublished)
tuning results by Douady and Hubbard
(compare \cite{mi2}).\vspace{1,5mm}

Whereas symmetry in the dynamic plane can mostly be explained 
by the action of the quadratic map, the situation is more complicated
in the parameter space. Often there is a correspondence between 
local structure in the dynamic plane and in the parameter space
which helps to understand a special symmetry in the parameter space.

Typical examples are the local similarities about 
Misiurewicz points found by Tan Lei
(see \cite{tl}): The neighborhoods of a Misiurewicz point $c$
in the Mandelbrot set and in the corresponding Julia set are
asymptotically similar in the Hausdorff metric.

Roughly speaking, this provides infinitely
many points in $M$ with a `local rotation symmetry'.\vspace{1,5mm}

The present paper deals with symmetries whose nature is a combinatorial
one. In particular, it proves two statements conjectured
by Lau and Schleicher (see \cite{laus2,schlei}).
The first of the conjectures, the {\em Correspondence Principle},
relates combinatorial structure in the dynamic plane and the parameter space
and, as in the above example, it forces a symmetry in the parameter space.
The latter symmetry has an exact description by the
second conjecture of Lau and Schleicher,
the {\em Translation Principle}.\vspace{1,5mm}

To put us in a position to formulate the two principles, it is necessary
to recall some quadratic iteration theory.
Beside the standard paper \cite{dh1} by Douady and Hubbard, 
the references are
\cite{bea,bra,carl,la,mi1,st,mi5,schlei2}. 
\paragraph{Formulation of the Translation Principle.}
For a given complex parameter $c$, 
let $p_c$ denote the quadratic map acting on the complex plane
by $p_c(z)=z^2+c$.
The {\em filled-in Julia set} $K_c$ of
$p_c$ is defined to consist of all points with bounded
orbits, and its boundary $J_c=\partial K_c$ in the plane
is said to be the {\em Julia set} of $p_c$.

The {\em Mandelbrot set} $M$ is the set of all complex 
$c$ with connected Julia set.
It contains the set $H$ of all $c$ for which $p_c$ possesses
an {\em attractive} periodic orbit of some period $m$:
The {\em multiplier} of the orbit, i.e.\ the derivative
of the first return map $p_c^m$ at any point of the orbit,
has absolute value less than 1. If such an orbit
exists for $p_c$, it is unique.\vspace{1,5mm}

A connectedness component $W$ of $H$
is called {\em hyperbolic component} of the Mandelbrot
set. Its period ${\sc Per}(W)$ is defined to be the
period of the unique attractive periodic orbit
for $p_c; c\in W$, which indeed does not depend on $c\in W$. 
The only hyperbolic component of period 1 is called the
{\em main hyperbolic component}.

The {\em multiplier map} which assigns to each point $c$
of a given hyperbolic component $W$ the multiplier of the
attractive orbit for $p_c$
forms a conformal isomorphism
from $W$ onto the unit disk and extents continuously
to a homeomorphism from the closure of $W$ onto the 
closed unit disk. The point
in $\partial W$ which is mapped to $e^{2\pi\nu i}$
is said to have {\em internal} angle $\nu\in[0,1[$.
Some special points in the closure of a hyperbolic component $W$
of period $m$
play a crucial role in understanding the combinatorial structure of $M$:
\begin{enumerate}
\item[P1)] The {\em center} $c_W$ of $W$ is the unique point  
mapped to $0$ by the multiplier map.
\item[P2)] The {\em root} $r_W$ of $W$ is the point in $\partial W$
with internal angle 0. If $m>1$, at $r_W$ the Mandelbrot set splits:
$M\setminus \{r_W\}$ consists of
two connectedness components.
\item[P3)] A {\em bifurcation point} of $W$ is a point in $\partial W$
with internal angle $\frac{p}{q}$ ($p$ and $q$ are required to
be relatively prime). Such a point and only such one is the root $r_U$
of a hyperbolic component $U$ bifurcating from $W$:
${\sc Per}(U)=qm$, and the component of $M\setminus \{r_U\}$
containing $U$ is called the $\frac{p}{q}$-{\em sublimb} of $W$.
\end{enumerate}
\begin{defi}\label{orderp}
{\rm (Order and Visibility in parameter space)}\\
Let $W,U$ be hyperbolic components. 
Then $W\prec U$ if there exists a simple curve 
intersecting 
$M$ in exactly one point and dividing the complex 
plane into two open parts,
one containing $U$ and the other one $W$ and the 
main component.
 
$U$ is said to be {\em visible} from $W$
if $W\prec U$ and there is no hyperbolic
component $V$ of period less than the period of $U$ with
$W\prec V\prec U$.
\end{defi}

We want to formulate the Translation Principle now:
In the trees we consider, the
`nodes' are hyperbolic components labeled by their periods
and the `edges' are defined by the successor relation
induced by $\prec$. 
Equivalence of such trees
includes the labeling.\vspace{3mm}\\
{\bf Translation Principle:}
{\it Let $W$ be a hyperbolic component of period $m$.
Then the trees of hyperbolic components visible from $W$
in any two sublimbs of denominators $q_1$ and $q_2$ coincide,
including the embedding into the plane, when all periods of
visible components in the first sublimb are increased by 
$(q_2-q_1)m$.}\vspace{3mm}

\begin{figure}
\begin{minipage}{105mm}
\epsfxsize 105mm
\epsffile{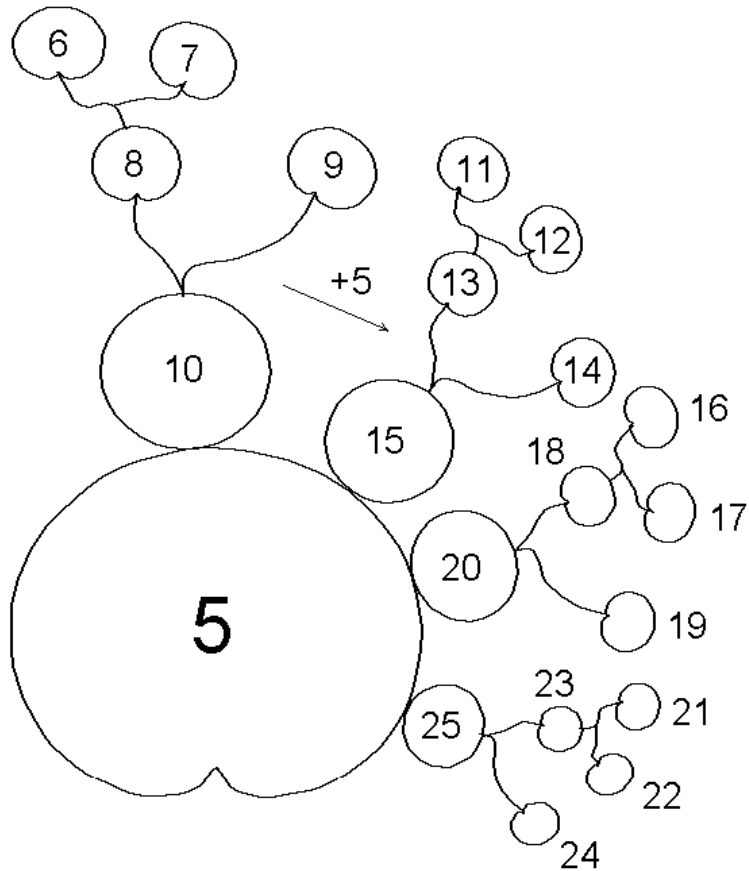}
\end{minipage}
\hfill
\begin{minipage}{47mm}\qquad\vspace{9cm}
\caption{The Trans\-la\-tion Prin\-ciple}\label{tpfig}
\end{minipage}
\end{figure}

Figure \ref{tpfig} demonstrates the Translation Principle. A 
hyperbolic component of period 5 together with the trees of the visible  
hyperbolic
components in sublimbs of denominators $q=2,3,4,5$ are shown.
\paragraph{The Correspondence Principle - a bridge between
dynamic plane and parameter space.}
The Translation Principle is forced by a
similar statement in the dynamic plane which is
caused by the action of a quadratic map,
and Lau and Schleicher's Correspondence Principle
translates between dynamic objects and objects 
in the parameter space.\vspace{1,5mm}

Fix a hyperbolic component $W$ of period $m$ and its 
center $c=c_W$. The quadratic map linking dynamic structure
and structure in the parameter space is $p_c$,
for which the {\em critical point} 0
and the {\em critical value} $c$ lie on a common periodic
orbit of period $m$. (We have the {\em superattractive case}, i.e.\
the multiplier of this orbit is 0.) 

The interior of $K_c$ decomposes into countably many
connectedness components, the {\em bounded Fatou components}.
Subsequently, `Fatou components' should always read 
`boun\-ded Fatou components'.
There are the {\em critical (Fatou) component}
and the {\em critical value
(Fatou) component} $A$ containing the
critical point and the critical value, respectively.
$A$ is the dynamic object assigned to $W$.\vspace{1,5mm} 

More general, in the dynamic plane the Fatou components 
take the role of the hyperbolic components.
The set of them is invariant under the action of $p_c$.
Moreover, $A$ has period $m$, and
for each Fatou component $F$ there exists an $n$
with $p_c^n(F)=A$.
The minimal $n$ with this property is called the {\sc Step} of
$F$. (Lau and Schleicher \cite{laus2} defined the {\sc Step}
of a precritical point instead of a Fatou
component in an analogous way. This applies to
more general cases,
but in the context of our paper there is no substantial difference because
each Fatou component contains exactly one backward 
iterate of $c$.)\vspace{1,5mm}
 
We want to complete the analogy between $W$ and $A$.
There is a unique homeomorphism from the closure of 
$A$ onto the closed unit disk which is conformal on $A$
and conjugates $p_c^m$ to the usual quadratic map
$p_0$. The unique point mapped to
$e^{2\pi\eta i}$ is said to have the 
{\em internal angle} $\eta\in [0,1[$.
\begin{enumerate}
\item[D1)] $c$ can be thought as the {\em center} of $A$. 
\item[D2)] The {\em dynamic root} $r_A$ of $A$ is the point in $\partial A$
with internal angle 0. At $r_A$ the filled-in Julia set splits:
$K_c\setminus \{r_A\}$ consists of finitely many components.
\item[D3)] A {\em bifurcation point} $x$ of $A$ is a point in $\partial A$
with internal angle $\frac{a}{2^{q-1}}$ 
for $q\geq 2$ and odd $a$:
The point $x$
is (first) mapped to $r_A$ after $q-1$ iterates of $p_c^m$,
the set $K_c\setminus \{x\}$ consists of finitely many components,
and the union of those which do not contain $A$
is called a $\frac{a}{2^{q-1}}$-{\em sublimb} of $A$
(compare \cite{laus2,schlei}).
\end{enumerate}
\begin{defi}\label{orderd}
{\rm (Order and Visibility in dynamic plane)}\\
Let $Q,R$ be Fatou components. 
Then $Q\prec R$ 
if there exists a simple curve intersecting 
$K_c$ in exactly one point and dividing 
the complex plane into two open parts,
one containing $R$ and the other one $Q$ and the 
critical component. 
 
$R$ is said to be {\em visible} from $W$
if $Q\prec R$ and there is no Fatou
component $S$ of {\sc Step} less than the {\sc Step} of $R$ with
$Q\prec S\prec R$.
\end{defi}

Let us formulate the Correspondence Principle.
The trees considered for the dynamic plane are
analogous to those for the parameter space,
and they are labeled by the {\sc Step}.
However, now 
we have trees in a generalized sense.
They need not have a root.
\vspace{3mm}\\
{\bf Correspondence Principle:}
{\it Let $W$ be a hyperbolic component of {\em period} $m$, 
let $c$ be its center
and $A$ be the critical value component of $p_c$.
Then, for a given sublimb $Subl_W$ of $W$ and of denominator $q$
and a given sublimb $Subl_A$ of $A$ and of denominator $2^{q-1}$
the following trees coincide,
including the labeling and the embedding in the plane:
\begin{enumerate}
\item[1.] The tree of all hyperbolic components in $Subl_W$
visible from $W$ and different from the hyperbolic component bifurcating
directly from $W$,
\item[2.] the tree of all
Fatou components visible from $A$ in $Subl_A$
of {\sc Step} less than $mq$.
\end{enumerate}}
\section{Thinking in laminations}
The nature of the Translation and Correspondence Principles
is a combinatorial one,
and there are different ways to extract
the combinatorial part from quadratic dynamics.
We want to use Thurston's concept of a lamination
(see \cite{th}).\vspace{1,5mm}

By $T$ we denote the unit circle, which we identify with the interval
$[0,1[\hspace{1mm}($`$=$'$R/Z)$ via $\beta\longleftrightarrow e^{2\pi\beta i}$.
Further, $h$ denotes the angle-doubling map 
$\beta\in T\longrightarrow 2\beta\bmod 1$ and $\beta_1\beta_2$
the chord with given ends $\beta_1,\beta_2$ in the unit circle $T$.
By a chord we understand a straight 
line contained in the closed unit disk with 
ends in the unit circle,
but without change of any statement below, a chord
can alternately be considered as the Euclidean closure
of a hyperbolic geodesic in the open unit disk.
(For our illustrations we use `hyperbolic' chords.)

By the {\em length} of a chord $B=\beta_1\beta_2$ and the {\em distance} of 
$\beta_1$ and $\beta_2$ we understand 
the minimum length of the two arcs in $T$ connecting
$\beta_1$ and $\beta_2$, where the whole circle has length 1.
The action of $h$ is extented to chords $B=\beta_1\beta_2$
by $h(B)=h(\beta_1)h(\beta_2)$.\vspace{1,5mm}

A chord $B=\beta_1\beta_2$ shorter than $\frac{1}{2}$
divides the unit disc into a smaller and a bigger open part.
$\beta\in T$ is said to be {\em between} $\beta_1$
and $\beta_2$ if it lies in the smaller part, and a subset of the unit disc 
is said to be {\em behind} 
$B$ if at least one of its points lies in the smaller,
but no one in the bigger part of the disk.
Similarly, a set is defined to be {\it between} two
disjoint chords $B_1,B_2$ 
if it lies in the closed part of the disk bounded by $B_1$ and $B_2$
but at least one of its points does not belong 
to $B_1$ or to $B_2$.

Finally, we shall use the following diction: Two chords $S_1,S_2$ {\it cross}
each other if they are different, but have a common interior point.
They are {\it separated} 
by a third chord $S_3$ if $S_3$ is between $S_1$ and $S_2$
but $S_2$ not between $S_1$ and $S_3$.
\begin{table}[h]
\begin{center}
\begin{tabular*}{14,8cm}{|p{6,5cm}@{\hspace{1,5mm}\extracolsep\fill}|p{7,5cm}|}\hline
complex plane & lamination model\\ \hline\hline
hyperbolic component $W$ of period $m>1$ & 
leaf $S$ of period $m$ in ${\cal B}_\ast$ \\ \hline
Fatou component of {\sc Step} $n$ & dynamic pair of {\sc Step} $n$\\ \hline
sublimb of denominator $q$ in para\-me\-ter space, 
therein: main antenna tip 
of the tuned Mandelbrot set
corresponding to $W$
&
region behind immediately visible
$B\in{\cal B}_\ast$ of period 
$qm$, preperiodic leaf $R_B$\\ \hline
sublimb of denominator $2^{q-1}$ in dynamic plane&
region behind 
$L_{{\bf v}s_1{\bf v}s_2\ldots {\bf v}s_{q-2}{\bf v}e}(S)$
for 0-1-word $s_1s_2\ldots s_{q-2}$
\\ \hline
\end{tabular*}
\end{center}
\centerline{Dictionary between complex plane and laminations}
\end{table}

In the following, a {\it lamination} is a set of non-degenerate
chords (i.e. chords with different ends)
which do not cross each other. (Here the concept of a lamination is used
in a more general way than Thurston's original one.) 
The elements of a lamination are called its {\em leaves},
and in the concrete situation
their ends correspond to neighboring external rays landing at
the same point of
the Mandelbrot set or Julia set.
By Lindel\"of's theorem the curves in the Definitions
\ref{orderp} and \ref{orderd} can be considered
as the union of two external rays and their common
landing point (but the concept of an external ray was
not necessary to formulate the principles.)

The above dictionary translating objects in the complex plane
into the correspon\-ding objects in the lamination setting should be an
orientation for the reader.
\unitlength1cm
\begin{figure}[h]
\begin{minipage}{7cm}
\begin{picture}(7,7)
\epsfxsize 70mm
\epsffile{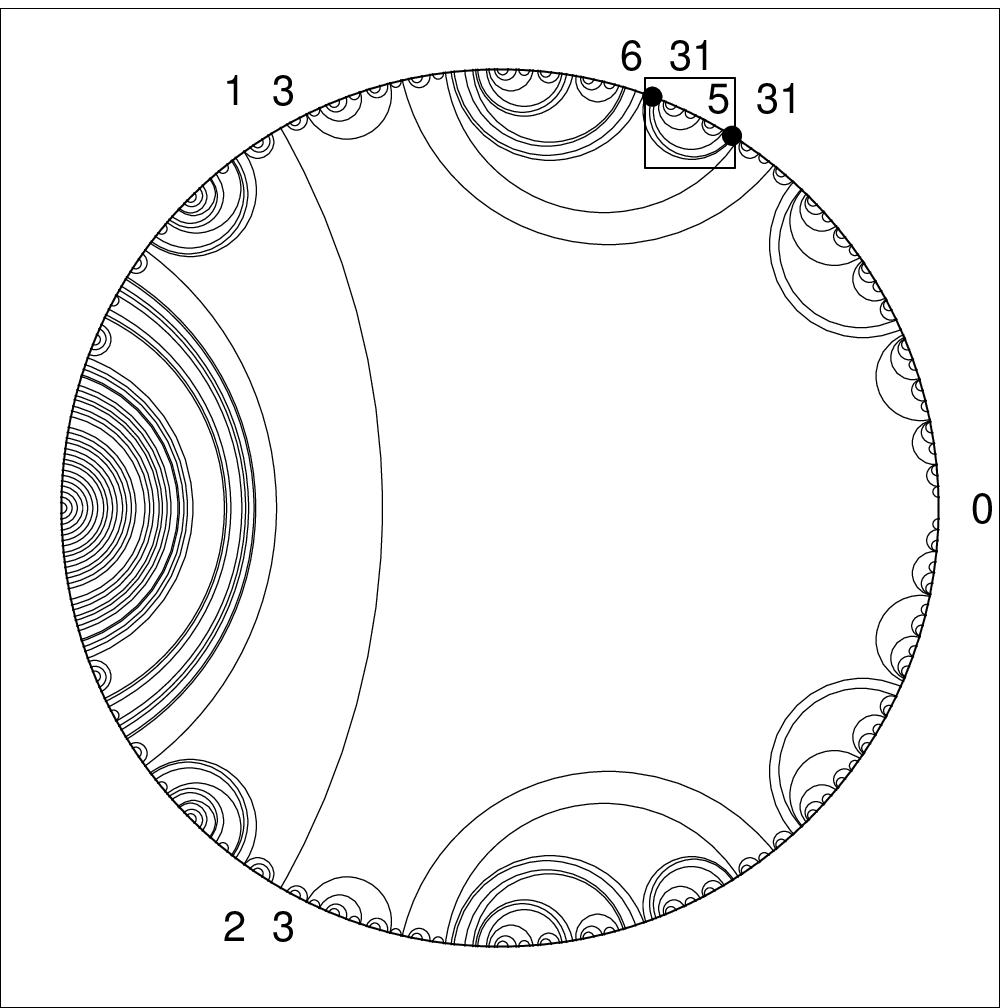}
\end{picture}
\caption{The periodic parameter lami\-nation ${\cal B}_\ast$}\label{bast}
\end{minipage}
\hfill
\begin{minipage}{7cm}
\begin{picture}(7,7)
\epsfxsize 70mm
\epsffile{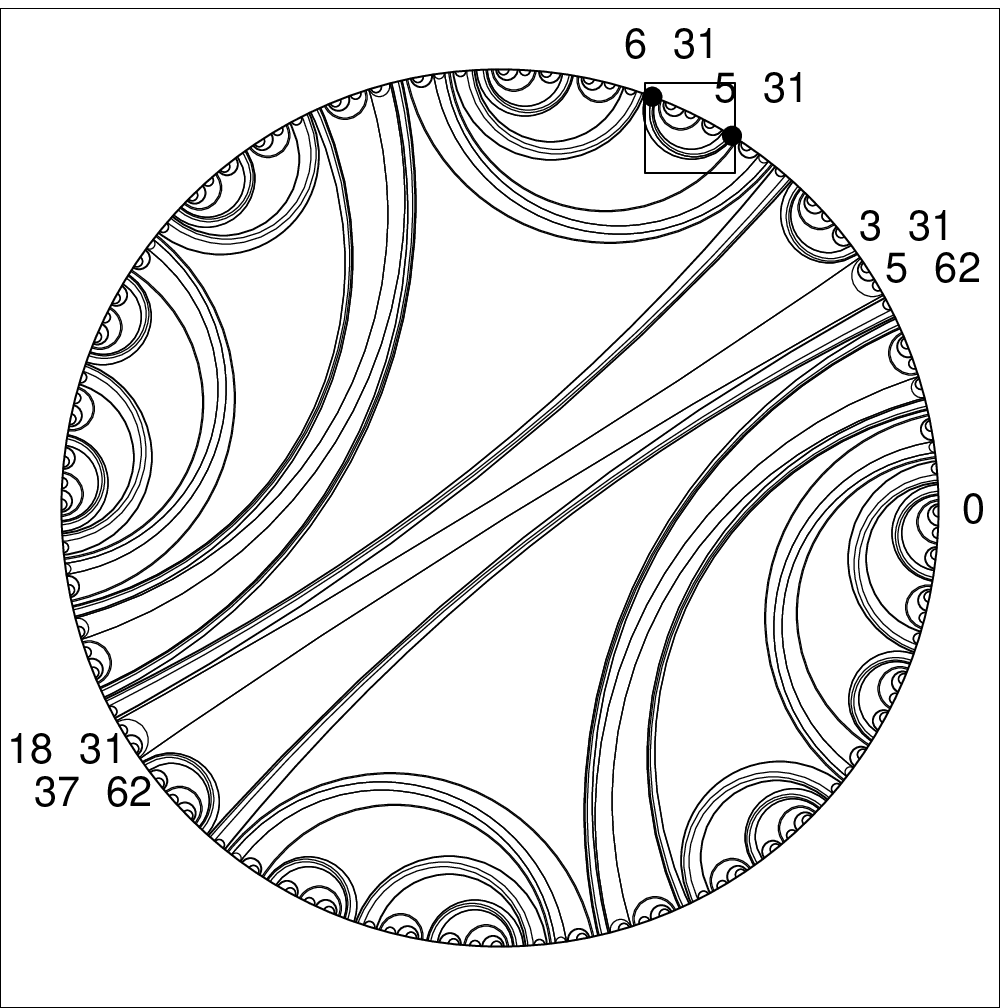}
\end{picture}
\caption{The generating 
dynamic la\-mi\-nation ${\cal B}_\ast(S)$
for $S=\frac{5}{31}\frac{6}{31}$}\label{dbast}
\end{minipage}
\end{figure}\\
{\bf a)} The {\bf periodic parameter lamination ${\cal B}_\ast$}
(see Figure \ref{bast})
consists of all {\it periodic parameter leaves} of some period $m>1$:
Their two ends are periodic of period $m$ (with respect to $h$) 
and form the external angles of the root of a hyperbolic component 
having period $m$. (By Douady and Hubbard \cite{dh1}, 
such roots have exactly two periodic external angles
and each periodic angle is an external angle of such a root.
For this also compare \cite{mi5,schlei2}.)

So each leaf in ${\cal B}_\ast$ of some period describes
a non-main hyperbolic component of the same period, and vice versa.
Let us express visibility of hyperbolic components
within the lamination concept. (`immediately visible' indicates
the direct bifurcation of a hyperbolic component.)
\begin{defi}\label{visietc} {\rm (Parameter visibility)}\\
Let $S\in {\cal B}_\ast$.
Then
\begin{enumerate}
\item[(i)]
$B\in {\cal B}_\ast$ 
is said to be {\em immediately visible} (from $S$) if 
$B$ is behind $S$, and
$S,B$ are not separated 
by a leaf of
${\cal B}_\ast$. 
\item[(ii)]
$B\in {\cal B}_\ast$ 
is said to be {\em visible} (from $S$) if 
$B$ is behind $S$, and
$S,B$ are not separated 
by a leaf of
${\cal B}_\ast$ whose period is less than the period
of $B$.
\end{enumerate}
\end{defi}

The meaning of `visible'
is not changed when
one substitutes the word `less' by `less or equal' in (ii).
This follows from a well known result of Lavaurs
\cite{la} (compare also Section 5 in \cite{bk2}):\vspace{3mm}\\
{\bf Lavaurs' Lemma} {\it Two leaves in
${\cal B}_\ast$ of the same period
are separated by a leaf of lower period.}\vspace{3mm}

The {\bf entire parameter lamination
${\cal B}$} consists of all {\em parameter leaves}
defined to be an element or a (non-degenerate) 
accumulation chord of ${\cal B}_\ast$.
By an accumulation chord of some set of chords ${\cal A}$
we understand the limit of a converging sequence of mutually different
elements taken from ${\cal A}$.
Note that ${\cal B}$ together with further one-point leaves
is just Thurstons `quadratic minor lamination'
(see \cite{th}).
By (9) in \cite{bk2}, we have the simple 
geometric characterization of the elements of ${\cal B}$
following now.
(What in the present paper is denoted ${\cal B}_\ast$ respectively 
${\cal B}$,
is ${\cal S}_\ast$ respectively ${\cal S}$ 
in \cite{bk2}.)\vspace{1,5mm}
\begin{eqnarray}\label{bdefi}
S=\alpha\gamma\in\cal B&\Longleftrightarrow&
\mbox{The iterates of }S
\mbox{ do not cross each other, and they}\nonumber\\
&&\mbox{do not cross the diameters }\frac{\alpha}{2}\frac{\alpha+1}{2}\mbox{ and }
\frac{\gamma}{2}\frac{\gamma+1}{2}. 
\end{eqnarray}

{\bf b)}
Fix a hyperbolic component $W$ of period $m>1$ with center $c$.
It is well known that the periodic parameter leaf
$S=\alpha\gamma\in{\cal B}_\ast$ corresponding to the root of $W$
generates the lamination belonging to $J_c$:

$S$ has a unique preimage with periodic ends denoted by $\dot{S}$
and being equal to $h^{m-1}(S)$ and one with preperiodic
ends denoted by $\ddot{S}$. The
leaves $\dot{S},\ddot{S}$ are the 
longest in the lamination, and it holds
$\{\dot{S},\ddot{S}\}=
\{\frac{\alpha}{2}\frac{\gamma+1}{2},
\frac{\alpha+1}{2}\frac{\gamma}{2}\}$.
The {\bf generating dynamic lamination ${\cal B}_\ast(S)$}
is defined to be the set of all chords $R$ such that
some iterate of $R$ is equal to $S$, no iterate of $R$ is longer
than the leaves $\dot{S}$ and $\ddot{S}$ 
or crosses one of them, 
and if an iterate of $R$ lies between the leaves $\dot{S}$ and $\ddot{S}$,
it has no end in common with them.

To $J_c$ there corresponds the {\bf entire dynamic 
lamination ${\cal B}(S)$} defined to be 
the set ${\cal B}_\ast(S)$ together with its non-degenerate
accumulation chords.
It is important to note that no chord in ${\cal B}(S)$
between $\dot{S}$ and $\ddot{S}$ has a point in common 
with $\dot{S}$ and $\ddot{S}$
and that no iterate of $S$ lies
between $\dot{S}$ and $\ddot{S}$.\vspace{1,5mm}

Figure \ref{dbast} shows the generating 
dynamic la\-mi\-nation ${\cal B}_\ast(S)$
for the leaf $S=\frac{5}{31}\frac{6}{31}$, which has period
5. To $S$ there corresponds
the hyperbolic component of period 5 sketched 
in Figure \ref{tpfig}, and it holds
$\dot{S}=\frac{3}{31}\frac{18}{31}$ and $\ddot{S}=\frac{5}{62}\frac{37}{62}$. 
\paragraph{Some statements on the structure of ${\cal B}(S)$.}
We want to encode the leaves in the lamination
${\cal B}_\ast(S)$ belonging to $J_c$. 
For this we use {\em 0-1-words},
i.e.\ words of symbols taken
from the alphabet $\{0,1\}$, including the empty word.
As usual, and for any alphabet, the notations are as follows:
${\bf w_1w_2}$ for the concatenation of two words ${\bf w}_1,{\bf w}_2$,
further, ${\bf w}^n$ for the $n$-fold concatenation of a word ${\bf w}$,
and $\overline{\bf w}$ instead of ${\bf www}\ldots$.
Further, let $\alpha$ be the minimum of the ends of $S$
(as elements of $[0,1[$).\vspace{1,5mm}

For a given 0-1-word ${\bf w}=w_1w_2\ldots w_k$
define a map $L_{\bf w}$ on
${\cal B}(S)$ as follows:
$L_0(R)$ and $L_1(R)$ 
denote the unique preimages of $R\in {\cal B}(S)$
in ${\cal B}(S)$
with at least one end in the open interval
$]\frac{\alpha}{2},\frac{\alpha+1}{2}[$ and 
at least one end in the open interval
$]\frac{\alpha+1}{2},\frac{\alpha}{2}[$, respectively,
and $L_{\bf w}(R):=L_{w_1}(L_{w_1}\ldots (L_{w_k}(R)))$.
(Intervals are taken in counter-clockwise direction.)

Further, let $e^S$ be the unique symbol $e\in\{0,1\}$
with $\ddot{S}=L_e(S)$,
and let ${\bf v}^S$ be the unique 0-1-word ${\bf v}$
of length $m-1$ with $S=L_{\bf v}(\dot{S})$.
Then $\dot{S}=L_{1-e^S}(S)$ and $S=L_{{\bf v}^Se^S}(S)$.\vspace{1,5mm}

There is a one-to-one correspondence between the Fatou components
of $K_c$ and the infinite gaps of ${\cal B}(S)$:
The complement of $\bigcup {\cal B}(S)$
splits into connectedness components.
The closure of such a component is said to be an
{\em infinite gap} if it is bounded by infinitely many
leaves of ${\cal B}(S)$. In particular, there is the {\em critical
value gap} bounded by $S$ and lying behind $S$,
and corresponding to the critical value component.\vspace{1,5mm}

Now we provide two propositions which are basic for the following.
The first one presents a well known description 
of bifurcation from the critical value component, but in the language
of laminations (for a lamination proof see \cite{ke4}).
\begin{pro}\label{prebound}
{\rm (Preperiodic boundary leaves of the critical value gap)}\\ 
Let $S\in {\cal B}_\ast$ be of period $m$ and let
${\bf v}={\bf v}^S,e=e^S$.
Then all boundary leaves of 
the critical value gap
of ${\cal B}(S)$ different from $S$ are preperiodic. They coincide
with $L_{{\bf v}s_1{\bf v}s_2\ldots {\bf v}s_l{\bf v}e}(S)$
for the 0-1-words $s_1s_2\ldots s_l$ and are
first mapped to $S$ after $(l+1)m$ iterates.
If $l\geq 1$, then
the first $m-1$-th iterates 
of the whole closed interval behind
$L_{{\bf v}s_1{\bf v}s_2\ldots {\bf v}s_l{\bf v}e}(S)$ 
are disjoint to the whole open interval behind $S$, and 
$h^m$ maps the whole
closed interval behind 
$L_{{\bf v}s_1{\bf v}s_2\ldots {\bf v}s_l{\bf v}e}(S)$
to that behind 
$L_{{\bf v}s_2{\bf v}s_3\ldots {\bf v}s_l{\bf v}e}(S)$
homeomorphically.
\qquad $\Box$
\end{pro}

\begin{figure}
\begin{picture}(7,7)
\put(2.7,0.3){$L_{00101}(S)$}
\put(3.1,1.5){\line(0,-1){0.7}}
\put(3,2.1){8}
\put(2.5,3.7){6}
\put(5.5,3.1){7}
\put(6.4,3.2){9}
\epsfxsize 70mm
\epsffile{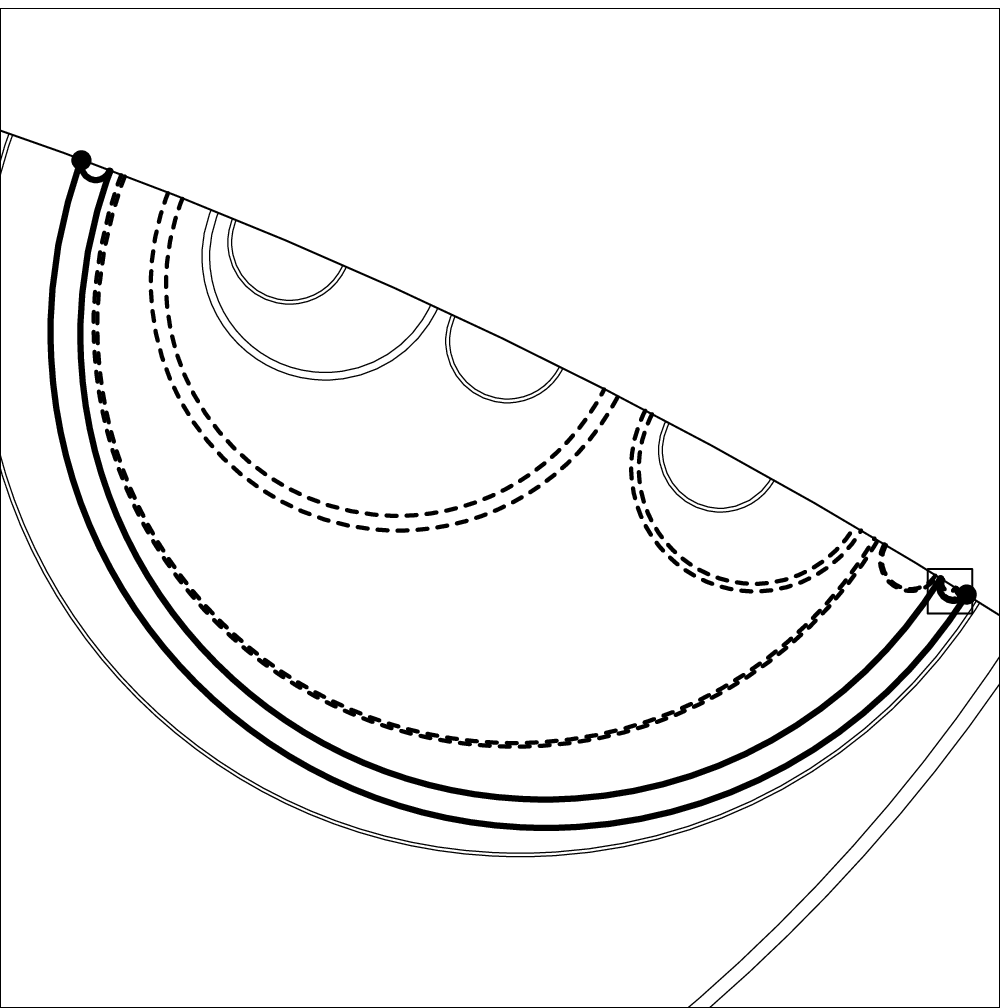}
\end{picture}
\hfill
\begin{picture}(7,7)
\put(1.8,0.9){$L_{0010100101}(S)$}
\put(2.8,2.4){\line(0,-1){1}}
\put(3.5,2.5){13}
\put(3.3,3.7){11}
\put(4.9,3.15){12}
\put(5.8,3.3){14}
\epsfxsize 70mm
\epsffile{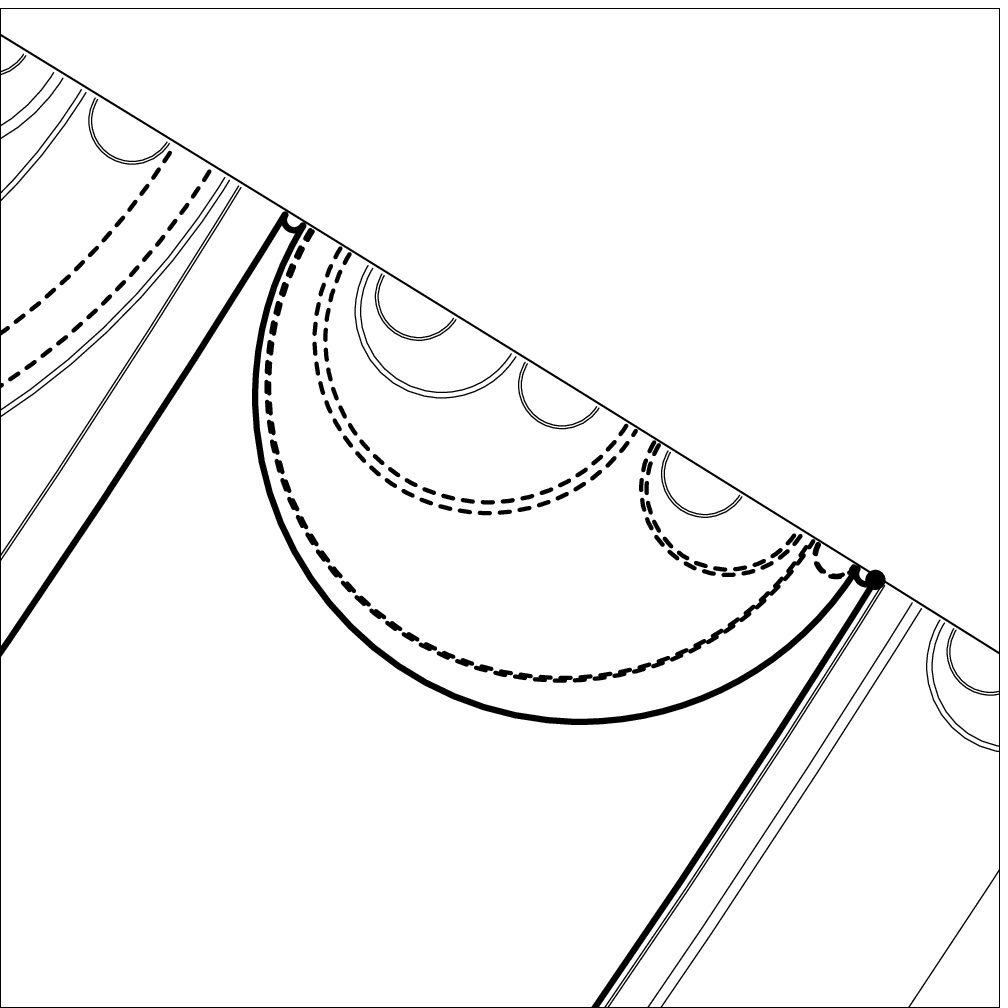}
\end{picture}
\caption{The Abstract Dynamical Translation Principle}\label{adtp}
\end{figure}

The left picture in Figure \ref{adtp} illustrates
${\cal B}_\ast(S)$ for $S=\frac{5}{31}\frac{6}{31}$
in a neighborhood of the critical value gap
(see window in Figure \ref{dbast}).
Boundary leaves of that gap are drawn as thick arcs. 
Note that ${\bf v}^S=0010$ and $e^S=1$.\vspace{1,5mm}

We come to the second proposition. Its first part
gives a well known statement encoding all Fatou components, 
in the language of laminations
(again see \cite{ke4}), 
and the second part
is not hard to show by using induction and the following simple statement
for $S\in {\cal B}_\ast$:
The leaves $\dot{S}$ and $\ddot{S}$ 
have length $a\geq\frac{1}{3}$ and span a `rectangle'
of `side'-lengths $a$ and $\frac{1}{2}-a$.
If some chord $Q$ has length $a'\leq a$,
then the length of $h(Q)$ is either $2a'$ or 
$1-2a'>\frac{1}{2}-a$.
\begin{pro}\label{longest}
{\rm (Infinite gaps in ${\cal B}(S)$)}\\
If $S\in {\cal B}_\ast$, then to each
0-1-word of length $n$ not ending with
${\bf v}={\bf v}^S$ there exists an infinite gap in ${\cal B}(S)$
different from the critical value one and vise versa.

The leaves $L_{{\bf w}0}(S)$ and
$L_{{\bf w}1}(S)$ 
form the longest `sides' of this gap:
If $S$ has length $d$, then
$L_{{\bf w}0}(S)$ and $L_{{\bf w}1}(S)$ 
are not shorter than
$\frac{d}{2^{n-1}}$ and
the two shortest `sides'
of the `rectangle' spanned by the ends of
$L_{{\bf w}0}(S)$ and $L_{{\bf w}1}(S)$ 
have length $\frac{d}{2^n}$.
\qquad $\Box$
\end{pro}

Proposition \ref{longest} shows that
the Fatou components (different from the critical value component)
can be described by two leaves of the corresponding lamination,
and this justifies the following definition:
\begin{defi}\label{dynpair} {\rm (Dynamic pairs and dynamic visibility)}\\
Let $S\in {\cal B}_\ast$ and 
$S_1,S_2\in {\cal B}_\ast(S)$.
Then $(S_1,S_2)$ is called a
{\em dynamic pair} if $S_1$ is not shorter than $S_2$
and there exists a 0-1-word ${\bf w}$ which does not
end with ${\bf v}^S$, such that
$\{S_1,S_2\}=\{L_{{\bf w}0}(S),L_{{\bf w}1}(S)\}$.     

By the {\sc Step} of a dynamic pair $(S_1,S_2)$ we understand
the minimal
number $n$ with $h^n(S_1)=h^n(S_2)=S$. (It exceeds the
length of the corresponding word ${\bf w}$ by one.)

A dynamic pair $(S_1,S_2)$ is said to be {\em visible} 
from a preperiodic boundary leaf $R$ of the critical value gap
if $S_2$ lies behind $R$ and there exists no dynamic pair $(Q_1,Q_2)$ 
whose {\rm Step} is less than the {\rm Step}
of $(S_1,S_2)$ such that $Q_1=R$ or $Q_1$ separates 
$R$ and $S_1$.
\end{defi}
{\bf Remark:} 
The statements in Proposition \ref{longest} concerning lengths of
chords show that
no dynamic pair visible from a preperiodic boundary leaf
lies between the members of a dynamic pair.
\paragraph{What shall we show substantially?}
Showing the Translation and Correspondence Principles, we
can assume that the hyperbolic component $W$ is different
from the main component. Namely, for the latter the statements 
are obvious: The filled-in Julia set for its center is the unit disk,
and it is easy to see that the only hyperbolic components 
visible from the main component, bifurcate 
directly from it (see also \cite{laus2}).\vspace{1,5mm}

We shall give our proofs within
the lamination models. 
(The necessary translations
are given in the table above.)
As mentioned at the beginning, the Translation Principle is forced
by a similar statement in the dynamic plane.
In the lamination setting, the latter is given 
by the following Theorem \ref{dyntheo} being
an immediate consequence of Proposition \ref{prebound}.
\begin{theo}\label{dyntheo} 
{\rm (Abstract Dynamical Translation Principle)}\\
Let $S\in {\cal B}_\ast$, 
let $m$ be the period of $S$ and
${\bf v}={\bf v}^S,e=e^S$.
Further, let $s_1s_2\ldots s_l$ be a non-empty 0-1-word.

Then $h^{lm}$ maps the whole closed interval behind
$L_{{\bf v}s_1{\bf v}s_2\ldots {\bf v}s_{l}{\bf v}e}(S)$
onto the whole closed interval behind $L_{{\bf v}e}(S)$
homeomorphically. In particular, $h^{lm}$ transforms  
the set of all dynamic pairs behind
$L_{{\bf v}s_1{\bf v}s_2\ldots {\bf v}s_{l}{\bf v}e}(S)$
into the set of all those behind $L_{{\bf v}e}(S)$,
decreasing the {\sc Step} by $lm$, and
a dynamic pair $(S_1,S_2)$ of {\sc Step} $n$
is visible from
$L_{{\bf v}s_1{\bf v}s_2\ldots {\bf v}s_l{\bf v}e}(S)$ iff
$(h^{lm}(S_1),h^{lm}(S_2))$ is visible from
$L_{{\bf v}e}(S)$
and is of {\sc Step} $n-lm$. \qquad $\Box$  
\end{theo}

Figure \ref{adtp} illustrates Theorem \ref{dyntheo}.
The magnification of the 
section in the window on the left side yields
the picture on the right side, 
and behind 
$L_{0010100101}(S)$
one sees the same
as behind $L_{00101}(S)$.
Namely, $h^5$ enlarges all behind the first leaf
by the factor $2^5$, in particular, the visible
dynamic pairs, indicated by 
thick dashed double-arcs and labeled by their {\sc Step}'s.\vspace{1,5mm}

According to the following proposition to each parameter
leaf $B\in {\cal B}_\ast$
visible from some $S\in {\cal B}_\ast$ there is assigned
a special (dynamic) boundary leaf $R_B$ of the critical value gap
of ${\cal B}(S)$. We will prove that proposition in Section \ref{proofs}.
For the reader familiar with the tuning construction
by Douady and Hubbard \cite{dh2} we refer to the characterization
of $R_B$ in the table above.
\begin{pro}\label{cobrb} 
{\rm (Correspondence between
immediately visible parameter leaves and special
boun\-dary leaves of the critical value gap)}\\
Let $S\in {\cal B}_\ast$ be of period $m$ and let $q\geq 2$.
Then for each $B\in {\cal B}_\ast$ immediately visible
from $S$ and of period $qm$ there exists a unique boundary leaf 
$R_B$ of the critical value gap 
satisfying
the following properties: $R_B$ is
first mapped to $S$ after $(q-1)m$ iterates, 
$R_B$ lies behind $B$, and
each leaf in ${\cal B}_\ast$ behind $B$ and visible
from $S$ lies behind $R_B$.
\end{pro}
Theorem \ref{dyntheo} shows that
the trees of Fatou components in two
sublimbs of the same deno\-mi\-nator coincide, 
including the labeling and the embedding in the plane, and
Proposition \ref{cobrb} establishes a relation
between the sublimbs in the parameter plane and
special sublimbs in the dynamic plane.
So the proof of the Correspondence Principle can be reduced
to showing the following
\begin{theo}\label{partheo} {\bf (Main Result)}
{\rm (Abstract Correspondence Principle)}\\
Let $S\in {\cal B}_\ast$ be of period $m$ and let
$B\in {\cal B}_\ast$ be immediately visible from $S$ and of
period $qm$ for $q\geq 2$. Then for each dynamic pair $(S_1,S_2)$
visible from $R_B$
of {\sc Step} less than
$qm$ there is a unique leaf $Q(S_1,S_2)\in {\cal B}_\ast$ visible
from $S$ and separating $S_1$ and $S_2$. Its period is equal to the
{\sc Step} of $(S_1,S_2)$.

Moreover, by $(S_1,S_2)\mapsto Q(S_1,S_2)$
the set of all those dynamic pairs is mapped bijectively
onto the set
of all leaves of ${\cal B}_\ast$ behind $B$
which are visible but not immediately visible
from $S$.
\end{theo}
\begin{figure}[h]
\begin{minipage}{70mm}
\begin{picture}(7,7)
\put(2.5,0.3){$S=\frac{5}{31}\frac{6}{31}$}
\put(0.1,0.8){$R_B=L_{00101}(S)$}
\put(2.5,5.8){$S_2=L_{001010}(S)$}
\put(2.5,2.45){$S_1=L_{001011}(S)$}
\put(5,5){$Q(S_1,S_2)$}
\put(0.7,6.4){$B$}
\put(4.4,4.5){\line(1,1){0.4}}
\put(1.3,4.4){\line(1,1){1.2}}
\put(3.2,1.3){\line(0,-1){0.6}}
\put(3.1,3.37){\line(0,-1){0.4}}
\put(1.7,2.25){\line(-1,-2){0.45}}
\put(0.8,6.0){\line(0,1){0.2}}
\epsfxsize 70mm
\epsffile{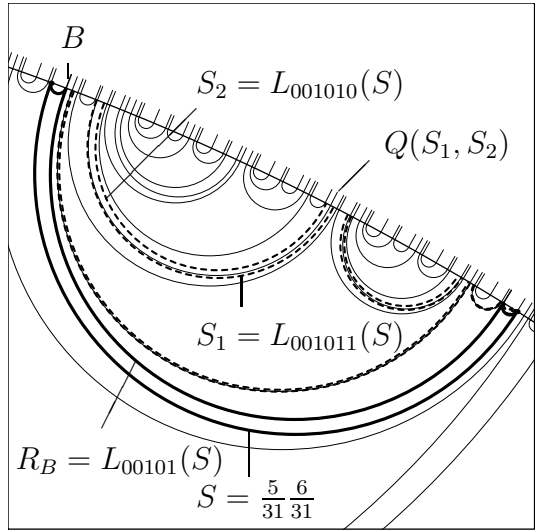}
\end{picture}
\end{minipage}
\hfill
\begin{minipage}{84mm}\qquad\vspace{5,7cm}
\caption{The Abstract Correspondence Principle}\label{acp}
\end{minipage}
\end{figure}

Figure \ref{acp} illustrates 
Proposition \ref{cobrb} and
Theorem \ref{partheo} for $S=\frac{5}{31}\frac{6}{31}$
and $q=2$.
Visible dynamic pairs 
are drawn in as in Figure \ref{adtp}. Beside them one sees
some elements of ${\cal B}_\ast$. To distinguish the latter
from the dynamic leaves, they are represented
by curves whose end points lie a little outside the unit disk.
$B$ has period $10$, and $R_B$ lies behind $B$ but
is very near to $B$.\vspace{1,5mm}

Obviously, the Translation Principle is an immediate consequence
of the Theorems \ref{dyntheo} and \ref{partheo}.
So the substantial part of the paper
is devoted to the proof of Theorem \ref{partheo}
(see Section \ref{proofs}).
We shall refer to it as the main result of the paper.
\section{Proofs}\label{proofs}
\paragraph{Kneading sequences.} 
For a given point $\alpha\in T$,
we call the sequence
\begin{eqnarray*}
s_1s_2s_3...
\quad\mbox{with}\quad s_i = \left\{\begin{array}{r@{\quad
\mbox{for}\quad }l} 0 & h^{i-1}(\alpha )\in\hspace{1,5mm}
]\frac{\alpha}{2},\frac{\alpha+1}{2}[\\ 
1 & h^{i-1}
(\alpha )\in\hspace{1,5mm}
]\frac{\alpha+1}{2},\frac{\alpha}{2}[\\ 
\ast & h^{i-1}(\alpha )\in
\{\frac{\alpha}{2},\frac{\alpha +1}{2}\} \end{array} \right.
\end{eqnarray*}
the {\em knea\-ding sequence} of 
$\alpha$. The concept of a kneading sequence has been introduced 
by Milnor and Thurston for the investigation of 
real quadratic maps (see \cite{mt}), and its use 
in the context of complex quadratic maps (in more or less explicit
form)
is due to different authors
(compare \cite{at,bk,bul,hus},
and note that the use of symbols in the kneading
sequence is not uniform in literature).\vspace{1,5mm}

The $n$-th symbol of the kneading sequence of $\alpha\in T$
is equal to $\ast$ iff $\alpha$ is
periodic and the period of $\alpha$ divides $n$.
An easy computation shows that this is satisfied iff
$\alpha$ has the form $\frac{a}{2^n-1}$ with $a\in\{0,1,\ldots 2^n-2\}$,
and exactly in the points having this form 
with $a\neq 0$, the $n$-th symbol
of the kneading sequence changes from 0 to 1.
(In particular, the period of $\alpha$ is the minimal $n$ for which
$\alpha$ has a representation of the form $\frac{a}{2^n-1}$.)
So the following is obvious (compare \cite{ke1,laus2}):
\begin{lem}\label{kslem}
Let $\beta_1,\beta_2$ be non-periodic points in $T$ 
with $\beta_1<\beta_2$. Then the $n$-th symbols of their
kneading sequences coincide iff the number of
points $\frac{a}{2^n-1}$ with $a=1,2,\ldots, 2^n-2$ in
$]\beta_1,\beta_2[$ is even.
\qquad $\Box$
\end{lem}

Fix some leaf $B=\alpha_1\alpha_2\in {\cal B}$. 
The number of points $\frac{a}{2^j-1}$ between 
$\alpha_1$ and $\alpha_2$ is even for each $j$.
This follows since the leaves in ${\cal B}_\ast$ are mutually disjoint
and since each $B\in {\cal B}$ with one periodic end 
belongs to ${\cal B}_\ast$
(e.g., see \cite{bk2}, Lemma 3 and Theorem 2(a),
also \cite{dh1,mi5,schlei2}).

Thus, if a sequence $(\alpha_i)_{i=1}^\infty$ 
of points behind $B$ converges to an end of $B$, 
then also the corresponding kneading sequences converge, and
the limit does not depend on
$(\alpha_i)_{i=1}^\infty$ and the end of $B$.      
In particular, if $B\in {\cal B}_\ast$ has period $m$,
then the limit is
$\overline{{\bf v}^Be^B}$, which
can easily be deduced since 
the $m-1$-th iterate of
$\alpha_i$ lies between $\dot{B}$ and $\ddot{B}$
for sufficiently large $i$. This justifies the following definition
(compare \cite{laus2}, Definition 5.3):
\begin{defi}\label{fpdef}
The {\em kneading sequence} of a leaf $B\in {\cal B}$ 
is defined to be that of its ends.
Moreover, if $B\in{\cal B}_\ast$ then 
$\overline{{\bf v}^B(1-e^B)}$ is called
the {\em kneading sequence just before} $B$
and $\overline{{\bf v}^Be^B}$ the {\em kneading sequence just behind} $B$.
If $B\in {\cal B}\setminus {\cal B}_\ast$, then the 
{\em kneading sequences just before} and {\em just behind} $B$
are the kneading sequence of $B$.
\end{defi}

If $B\in {\cal B}$ 
is separated from some $Q\in {\cal B}$ behind $B$
by a leaf 
of period $m$, and if $m$ is assumed to be
minimal with this condition, then 
the corresponding leaf is unique by Lavaurs' Lemma.
This shows that the number of periodic points of period $m$
in each of the two intervals between $B$ and $Q$
having one end in common with both $B$ and $Q$ is odd.
So by Lemma \ref{kslem}
one easily obtains
the following statement (compare \cite{laus2}, Proposition 5.4):
\begin{lem}\label{futpas}
{\rm (Symbolic description of parameter visibility)}\\
Let $Q\in {\cal B}$ be a leaf behind some $B\in {\cal B}$.
Then $B,Q\in {\cal B}$ are not separated by
a leaf in ${\cal B}_\ast$ of period less than or equal to $n$
iff the initial subwords of length $n$ of the kneading sequences
just behind $B$ and just before $Q$ coincide.
In particular, a leaf $Q\in {\cal B}_\ast$ of period $m$
is visible from $B\in {\cal B}_\ast$ iff 
${\bf v}^B$ is equal to the initial subword of length $m-1$ 
just behind $Q$.
It is immediately visible if furthermore the period of $B$ divides $m$.
\qquad $\Box$
\end{lem}

If $S\in {\cal B}_\ast$, then the elements of the 
forward orbit of a given leaf $R=\alpha_1\alpha_2\in {\cal B}(S)$ 
behind $S$ do not 
cross $\dot{S}$ and $\ddot{S}$ and are not longer than
$\dot{S},\ddot{S}$.
Thus $h^i(R)$ for $i=0,1,2,\ldots$ does not cross
$\frac{\alpha_1}{2}\frac{\alpha_1+1}{2}$ and
$\frac{\alpha_2}{2}\frac{\alpha_2+1}{2}$ iff
$h^{i+1}(R)$ does not separate $S$ and $R$.
With (\ref{bdefi}) and
statement (b) in Theorem 2 of \cite{bk2} (compare also
\cite{mi5} and \cite{ke4}) we have the following
\begin{pro}\label{firstcorres}
{\rm (Important dynamic chords lie in ${\cal B}_\ast$)}
For $S\in{\cal B}_\ast$
the following statements are valid:
\begin{enumerate}
\item[(i)]
Each boundary leaf of the critical value gap of ${\cal B}(S)$
and each leaf belonging to a dynamic pair visible from a
boundary leaf of the critical value gap 
lies in ${\cal B}$.
\item[(ii)] If $B\in {\cal B}_\ast$ separates
0 and $S$, then $B,\dot{B}$ and $\ddot{B}$ are contained
in ${\cal B}(S)$. \qquad $\Box$
\end{enumerate}
\end{pro}
\paragraph{Proof of Proposition \ref{cobrb}.} 
Let $S\in {\cal B}_\ast$ of period $m$ and of length $d$
with minimum end $\alpha$ be given.
The boundary of the critical value gap of ${\cal B}(S)$
consists of the boundary leaves and
further elements of $T$.
Contracting each boundary leaf 
to a point, one gets a space 
homeomorphic to $T$. 

There is a special map 
from $T$ onto that space, which preserves orientation and
conjugates $h$ and $h^m$. 
That map, which is denoted by $f^\alpha$ here and
was described in \cite{bk2}, is due to (the combinatorial
description of tuning by)
Douady (see \cite{dou}).\vspace{1,5mm}

$f^\alpha$ is defined as follows (for equivalence
to Douady's definition, see Proposition 2 in \cite{ke2}):
For $\beta\in T$ with binary expansion
$.b_1b_2b_3\ldots$, i.e.\ with 
$\beta =\sum_{i=1}^\infty b_i 2^{-i}$, let
\begin{eqnarray}\label{fo2}
f^\alpha (\beta)
=\alpha+(2^m-1)\ d\ \sum_{i=1}^\infty b_i2^{-im}.
\end{eqnarray}

$f^\alpha$ assigns to each $\beta$ with
unique binary expansion a point different
from all contracted boundary leaves.
If $\beta\in T$ is first mapped to 0 after $l$ iterates,
then the two binary
expansions $.b_1b_2\ldots b_{l-1}1\overline{0}$ and
$.b_1b_2\ldots b_{l-1}0{\overline 1}$ lead to
the two endpoints of a boundary leaf $f^\alpha(\beta)$, 
which according to (\ref{fo2}) has length
\begin{eqnarray*}
(2^m-1)\ d\
\left (2^{-lm}-\sum_{i=2}^\infty 2^{-(l+i)m}\right )=
(2^m-1)\ d\ \left (\frac{1}{2^{lm}}-\frac{1}{2^{lm}(2^m-1)}\right )=
\frac{2^m-2}{2^{lm}}\ d
\end{eqnarray*}
in the case $l>0$.\vspace{1,5mm}

Now fix some $B\in {\cal B}_\ast$ of period $qm$
being immediately visible from $S$. As well known,
the ends of such $B$ lie on a common orbit.
Further, let  
$\gamma,\delta\in T$ be the unique points
which are mapped by $f^\alpha$ to the ends of $B$.
They must be periodic of period $q$ and must lie
on a common orbit.
Moreover, the smaller closed interval $I$ with ends $\gamma$ and $\delta$
is not shorter than $\frac{1}{2^q-1}$. 

>From this one easily deduces $T=\bigcup _{i=0}^{q-1}h^i(I)$,
which implies the existence of some $\beta\in T$ first mapped
to 0 after $l$ iterates, were $0<l\leq q-1$. 
Let $R_B:=f^\alpha(\beta)$.\vspace{1,5mm}

Clearly, $R_B$ lies behind the leaf $B$,
which has length
$\frac{(2^m-1)^2}{2^{qm}-1}d$ by a result of
Lau and Schleicher (see Proposition 2.7 in \cite{laus2}). 
This and the formula
above show that $l=q-1$. Since $h^{i-1}(R_B)$ lies between
$\dot{S}$ and $\ddot{S}$ iff $m$ divides $i$,
one obtains that $R_B$ is first mapped to $S$ after
$(q-1)m$ iterates and that 
the kneading sequence of $R_B$ starts with
${\bf v}s_1{\bf v}s_2\ldots s_{q-1}{\bf v}$, where
${\bf v}={\bf v}^S$ and $s_1,s_2,\ldots,s_{q-1}\in\{0,1\}$.\vspace{1,5mm}

Now let a leaf $Q\in {\cal B}_\ast$ behind $B$ be visible from $S$.
Then its period is less than $qm$.
By Lemma \ref{futpas} it cannot separate
$B$ and $R_B$, and 
by Proposition \ref{firstcorres}(i) it does not cross $R_B$.
The difference between the lengths of $B$ and $R_B$ is
\begin{eqnarray*}
\left( \frac{(2^m-1)^2}{2^{qm}-1}-\frac{2^m-2}{2^{(q-1)m}}\right) d=
\frac{2^{qm}+2^{2m}-2^{m+1}}{2^{qm}(2^{qm}-1)}\ d<\frac{2d}{2^{qm}-1}<
\frac{1}{2^{qm}-1}.
\end{eqnarray*}
Therefore $Q$ must lie behind $R_B$.
\paragraph{Proof of the main result.} We are going to prove 
Theorem \ref{partheo} now.
{\sf For the whole proof} 
fix a leaf $S\in {\cal B}_\ast$ of a given period $m$, and let
${\bf v}={\bf v}^S=v_1v_2\ldots v_{m-1}$ and $e=e^S$.
(For some possibly helpful illustrations, we refer to Figures
\ref{dbast}, \ref{adtp} and \ref{acp}, which show the case
$S=\frac{5}{31}\frac{6}{31}$ with period 5 and
${\bf v}^S=0010,e^S=1$.)
At this place also recall the well known fact that
the kneading sequence just behind
$S$ has period $m$ (e.g., see \cite{laus2}).

First of all, let us verify the following simple statement:
\begin{lem}\label{perlem}
If for a given $i$ the subword of $\overline{{\bf v}e}$
starting with the $i$-th symbol and having length $m-1$
is equal to ${\bf v}$, then
$i=jm+1$ for some nonnegative integer $j$.
\end{lem}
Proof: Under the assumptions of the Lemma, the $i+m-1$-th symbol must 
be equal to $e$. This can be shown by counting the symbols in a sufficiently
long initial subword of the periodic sequence $\overline{{\bf v}e}$.
So, if $i$ were different from all $jm+1$ for  
$j=0,1,2,\ldots$, then the period of $\overline{{\bf v}e}$ would be 
less than $m$, which is false.
\qquad $\Box$\vspace{3mm}

We shall use the following simple statement, which implies that
for given 0-1-words
${\bf u},{\bf v}$ and ${\bf w}$ and 
$s\in \{0,1\}$, 
the leaf
$L_{s{\bf u}}(S)$ separates $L_{s{\bf v}}(S)$ from $L_{s{\bf w}}(S)$ 
iff $L_{\bf u}(S)$ separates $L_{\bf v}(S)$ from $L_{\bf w}(S)$.
(One implication is an immediate consequence
of the following statement, and the other one can
be shown indirectly, also applying that statement.)\vspace{3mm}\\
{\bf Order-invariance Principle:} {\it If three (or more)
points are contained in
an open semi-circle, then the action of $h$ does not change their
circular order. \qquad $\Box$}\vspace{3mm}

Before starting the proof of Theorem \ref{partheo}, 
let us make an agree\-ment: Leaves $S_1,S_2,\ldots,S_k$
are said to be {\em in the arrangement} $S_1,S_2,\ldots,S_k$ if,
for each $i=3,4,\ldots, k$, the leaf $S_{i-1}$ separates 
$S_{i-2}$ and $S_i$, or
$S_{i-2}=S_{i-1}$ or $S_{i-1}=S_i$.\vspace{4mm}\\
{\sl I. Visible dynamic pairs determine (not immediately) visible 
periodic parameter leaves:}\vspace{2mm}\\
We start looking at the dynamic pairs visible from
the longest preperiodic boundary leaf of the critical value gap of 
${\cal B}(S)$.  
\begin{lem}\label{finlem}
Let $(S^0_1,S^0_2)$ be a dynamic pair visible from $L_{{\bf v}e}(S)$
and having {\sc Step} $jm+k\geq 2$
with $j\in\{0,1\}$ and $0\leq k<m$.

Then $k\neq 0$, 
the kneading sequence of $S^0_1$ is equal to
$({\bf v}e)^j v_1v_2\ldots v_k
\overline{{\bf v}(1-e)}$, 
and the kneading sequence of $S^0_2$ is equal to
$({\bf v}e)^j v_1v_2\ldots v_{k-1}(1-v_k)
\overline{{\bf v}(1-e)}$.
(The kneading sequences of the leaves $S^0_1$ and $S^0_2$ 
are well-defined since by Proposition \ref{firstcorres}(i)
they belong to ${\cal B}$.)
\end{lem}
Proof: The arguments of the proof will be 
based on studying the mutual positions
of the iterates of the leaves 
$S=L_{{\bf v}(1-e)}(S),L_{{\bf v}e}(S),S^0_1,S^0_2$, which
are in the arrangement $S,L_{{\bf v}e}(S),S^0_1,S^0_2$.
Further, we use the statement that for $\alpha$ behind $S$ the 
$i$-th symbol of the kneading sequence is $e$ if
$h^{i-1}(\alpha)$ lies behind $\ddot{S}$, and $1-e$ if
$h^{i-1}(\alpha)$ lies behind $\dot{S}$. This follows
from $\ddot{S}=L_e(S),\dot{S}=L_{1-e}(S)$ and the fact that
$\frac{\beta}{2}\frac{\beta+1}{2}$ separates
$\ddot{S}$ and $\dot{S}$.\vspace{1,5mm}

We start with the case that $j=0$. If then $0\leq
i<k-1$, the leaves $h^i(S^0_1),
h^i(S^0_2)$ cannot be separated from $h^i(S)$
by $L_e(S)=\ddot{S}$ or
$L_{1-e}(S)=\dot{S}$.
Otherwise, by the Order-invariance Principle, 
$S^0_1$ or $S^0_2$ 
would be separated from $S$ by $L_{v_1v_2\ldots v_ie}(S)$ or
$L_{v_1v_2\ldots v_i(1-e)}(S)$, which is a contradiction to
visibility of $(S^0_1,S^0_2)$.
In particular, by the Order-invariance Principle,
the $k-1$-th iterates of $S,S^0_1,S^0_2$ are in the arrangement
$h^{k-1}(S),h^{k-1}(S^0_1),h^{k-1}(S^0_2)$.

To see that 
$S^0_1$ and $S^0_2$ have kneading sequences
$v_1v_2\ldots v_k\overline{{\bf v}(1-e)}$ and 
$v_1v_2\ldots v_{k-1}(1-v_k)
\overline{{\bf v}(1-e)}$,
one only needs that
$\{h^{k-1}(S^0_1),h^{k-1}(S^0_2)\}=
\{\dot{S},\ddot{S}\}$.
(Recall that the iterates of $S$ do not lie
between $\dot{S}$ and $\ddot{S}$.)

$j=1$ and $k=0$ is impossible. Indeed, otherwise
by use of the same arguments as in the case $j=0$
one would obtain
$\{S^0_1,S^0_2\}=\{S,L_{{\bf v}e}(S)\}$
from $\{h^{m-1}(S^0_1),h^{m-1}(S^0_2)\}=
\{\dot{S},
\ddot{S}\}$, which contradicts
the assumption that 
$(S^0_1,S^0_2)$ is a (visible) dynamic pair.\vspace{1,5mm}

Finally, let us come to the case $j=1$ and $k>0$.
Again from the arguments in the case $j=0$, it follows
that the leaves $S,S^0_1$ and
$S^0_2$ have kneading sequences
with common initial subwords of length $m-1$. 

To show that the $m$-th symbols
of all kneading sequences considered are equal to $e$,
also take into consideration the leaf $L_{{\bf v}e}(S)$.
Using the Order-invariance Principle, one easily sees
that the leaves $\dot{S}=h^{m-1}(S),\ddot{S}=
h^{m-1}(L_{{\bf v}e}(S)),h^{m-1}(S^0_1),h^{m-1}(S^0_2)$
are in the arrangement 
$\dot{S},\ddot{S},
h^{m-1}(S^0_1),h^{m-1}(S^0_2)$,
which shows that both $S^0_1$ and $S^0_2$
have kneading sequences whose $m$-th symbols are equal to $e$.\vspace{1,5mm}

For the rest, one argues as at the beginning
of the proof: 
If $0\leq i<k-1$, then the leaves $h^{m+i}(S^0_1),
h^{m+i}(S^0_2)$ cannot be separated from $h^{m+i}(S)=h^i(S)$
by $L_e(S)=\ddot{S}$ or
$L_{1-e}(S)=\dot{S}$. 
Otherwise,  
$S^0_1$ or $S^0_2$ 
would be separated from $S$ by $L_{{\bf v}ev_1v_2\ldots v_ie}(S)$ or
$L_{{\bf v}ev_1v_2\ldots v_i(1-e)}(S)$, but 
$L_{{\bf v}ev_1v_2\ldots v_ie}(S)$ and by Lemma \ref{perlem} 
$L_{{\bf v}ev_1v_2\ldots v_i(1-e)}(S)$ 
belong to a dynamic pair having {\sc Step} less than $m+k$.

Moreover,
$\{h^{m+k-1}(S^0_1),h^{m+k-1}(S^0_2)\}=
\{\dot{S},
\ddot{S}\}$,
and the $m+k-1$-th iterates of $S,S^0_1,S^0_2$ are in the arrangement
$h^{m+k-1}(S),h^{m+k-1}(S^0_1),h^{m+k-1}(S^0_2)$.
\qquad $\Box$\vspace{3mm}\\
Now we are able to prove the first part of Theorem \ref{partheo}, 
which is covered by the following
\begin{cor}\label{maincor}
Let $B\in {\cal B}_\ast$ be immediately visible from $S$
and of period $qm$ for some $q\geq 2$.
Further, let $(S_1,S_2)$ be a dynamic pair 
in ${\cal B}_\ast(S)$ visible from $R_B$ and having {\sc Step} 
$jm+k\geq 2$ with $j<q$ and $0\leq k<m$. Then the following statements
are valid:
\begin{enumerate} 
\item[(i)] $j=q-2$ or $j=q-1$, and $k\neq 0$.
\item[(ii)] 
The kneading sequence of $S_1$ is equal to
$({\bf v}e)^j v_1v_2\ldots v_k
\overline{{\bf v}(1-e)}$, and
the kneading sequence of $S_2$ is equal to
$({\bf v}e)^j v_1v_2\ldots v_{k-1}
(1-v_k)\overline{{\bf v}(1-e)}$. 
\item[(iii)] There exists a unique leaf $Q(S_1,S_2)\in{\cal B}_\ast$
of period $jm+k$
between $S_1$ and $S_2$. 
It has kneading sequence
$(\overline{{\bf v}e)^jv_1v_2\ldots v_{k-1}\ast}$
and it is the only visible leaf in ${\cal B}_\ast$ separating
$S_1$ and $S_2$.
\end{enumerate} 
\end{cor}
Proof: 
$h^{(q-2)m}$ maps the leaf $R_B$ onto 
$L_{{\bf v}e}(S)$.
If $(S_1,S_2)$ is a dynamic pair
in ${\cal B}_\ast(S)$ visible from $R_B$
and having {\sc Step} $n=jm+k$,
then by Theorem \ref{dyntheo}
the dynamic pair $(S^0_1,S^0_2)$ with
$S^0_1=h^{(q-2)m}(S_1)$ and $S^0_2=h^{(q-2)m}(S_2)$
is visible from $L_{{\bf v}e}(S)$.
So one easily obtains (i) and (ii) from
Proposition \ref{prebound} and Lemma \ref{finlem}.\vspace{1,5mm}

(iii): Denote the two open intervals between $S_1$ and $S_2$
each having one end together with $S_1$ and one together
with $S_2$ by $I_1$ and $I_2$. By (ii) and Proposition \ref{longest},
these intervals are shorter than $\frac{1}{2^n}$. 
Thus, for each $l\leq n$, 
both $I_1$ and $I_2$ contain at most one periodic point of period $l$
(see above Lemma \ref{kslem}). 
By Lemma \ref{kslem}, there is one in $I_1$ and one in $I_2$ for $l=n$,
since in this case the kneading sequences of $S_1$ and $S_2$ have different
$n$-th symbols, but no one for $l<n$ since the $l$-th symbols 
coincide now.\vspace{1,5mm}

Let $Q=Q(S_1,S_2)$ be the leaf connecting the two periodic points
of period $n$ we obtained. $S_1,S_2\in {\cal B}$
forces that $Q\in {\cal B}_\ast$, and obviously the kneading sequence
of $Q$ is equal to $(\overline{{\bf v}e)^jv_1v_2\ldots v_{k-1}\ast}$.
Moreover, by Lemma \ref{futpas}
it is visible and, since $0<k<m$, not immediately visible from $S$.

What remains is to verify that there is no
other leaf in ${\cal B}_\ast$ visible from $S$ 
and separating $S_1$ and $S_2$. 
Assuming the opposite, the period of such a leaf 
must be greater than $n$ and so by visibility it must separate $S_1$
and $Q$.\vspace{1,5mm}

Let $r$ be the minimum of periods of all leaves in ${\cal B}_\ast$
separating $S_1$ and $Q$. Then by Lavaurs' Lemma 
there is a unique one of period $r$, and $r>n$.
Let ${\bf s}$ and ${\bf q}$ be the kneading sequence of $S_1$
and the kneading sequence just before $Q$, respectively.
Then by Lemma \ref{futpas}, the initial subwords of length
$r-1$ of ${\bf s},{\bf q}$ and $\overline{{\bf v}e}$
coincide, but the $r$-th symbols of ${\bf s}$ and ${\bf q}$
are different.

Thus from 
${\bf s}=({\bf v}e)^j v_1v_2\ldots v_k
\overline{{\bf v}(1-e)}$ and
${\bf q}=
(\overline{{\bf v}e)^jv_1v_2\ldots v_{k-1}v_k}$
one easily obtains $r\geq(j+1)m+k$, implying
$({\bf v}e)^jv_1v_2\ldots v_k{\bf v}=
({\bf v}e)^{j+1}v_1v_2\ldots v_k$. 
Therefore, we have
${\bf v}=v_{k+2}\ldots v_{m-1}ev_1\ldots v_k$, which
contradicts Lemma \ref{perlem}.
\qquad $\Box$\vspace{4mm}\\
{\sl II. All visible but not immediately visible periodic parameter
leaves are covered:}\vspace{2mm}\\
According to the remark below Definition \ref{dynpair}
the map $Q(\cdot)$ is injective. 
Now we show its surjectivity. The first step is given by
\begin{cor}\label{never}
Let $B\in {\cal B}_\ast$ be immediately visible from $S$.
Further, let $Q\in {\cal B}_\ast$ of period $n$ be 
visible from $S$ and behind $B$, and let $(S_1,S_2)$ be a dynamic pair
of {\sc Step} $k\leq n$ in ${\cal B}_\ast(S)$.
If $S_1$ separates $S$ and $Q$, then $k=n$
and $(S_1,S_2)$ is visible from $R_B$.
\end{cor}
Proof:
$Q$ lies behind $R_B$, and under the assumption that
$S_1$ separates $S$ and $Q$, the leaf $S_2$ must be behind $R_B$.
If $(S_1,S_2)$ were not visible from $R_B$, then one could
find a dynamic pair
$(S_3,S_4)$ visible from $R_B$
such that $S_1$ lies behind $S_3$, 
and by the remark below Definition \ref{dynpair}
also behind $S_4$.
So $Q(S_3,S_4)$ would be a visible leaf in ${\cal B}_\ast$ of period
less than $n$, which
contradicts visibility of $Q$.\qquad $\Box$\vspace{3mm}

Obviously, the proof of Theorem \ref{partheo} becomes complete 
by showing the following
\begin{lem}\label{finito}
Let $B\in {\cal B}_\ast$ be immediately visible from $S$.
Further, let $Q\in {\cal B}_\ast$ be visible from $S$ and behind $B$,
and let ${\bf w}={\bf v}^Q$.

Then $L_{{\bf w}0}(S)$ and $L_{{\bf w}1}(S)$
are separated by $Q$ and 
form a dynamic pair in ${\cal B}_\ast(S)$. 
\end{lem}
Proof: Let $n$ be the period of $Q$ and let
${\bf w}=w_1w_2\ldots w_{n-1}$. 
Then by Lemma \ref{futpas} and 
Lemma \ref{perlem} 
the word ${\bf w}$ does not end with ${\bf v}$,
thus the two leaves $L_{w_kw_{k+1}\ldots w_{n-1}}(\dot{S}),
L_{w_kw_{k+1}\ldots w_{n-1}}(\ddot{S})$ form a dynamic pair 
for all $k=1,2,\ldots ,n-1$. (The maps $L_{\bf u}$ for
0-1-words ${\bf u}$ are taken
with respect to $S$, i.e.\ all leaves considered belong to
${\cal B}_\ast(S)$.)\vspace{1,5mm}

According to Proposition \ref{firstcorres}(ii),
one has 
$S,\dot{S},\ddot{S}\in {\cal B}(Q)$, and
by induction one shows that
for all $l=1,2,\ldots ,n-1$
both $L_{w_lw_{l+1}\ldots w_{n-1}}(\dot{S})$ and
$L_{w_lw_{l+1}\ldots w_{n-1}}(\ddot{S})$ 
belong to ${\cal B}(Q)$.

The induction step is as follows:
Assume that $L_{w_{k+1}w_{k+2}\ldots w_{n-1}}(\dot{S}),
L_{w_{k+1}w_{k+2}\ldots w_{n-1}}(\ddot{S})$ are elements of ${\cal B}(Q)$
but at least one of the leaves $L_{w_kw_{k+1}\ldots w_{n-1}}(\dot{S})$
and $L_{w_kw_{k+1}\ldots w_{n-1}}(\ddot{S})$ does not belong to
${\cal B}(Q)$. Then the infinite gap in ${\cal B}(S)$ defined
by the latter leaves lies between $\dot{S}$
and $\ddot{S}$, but none of its boundary leaves separates
$\dot{S}$ and $\ddot{S}$. 

Denote the longer of the leaves $L_{w_kw_{k+1}\ldots w_{n-1}}(\dot{S})$
and $L_{w_kw_{k+1}\ldots w_{n-1}}(\ddot{S})$ by $R$.
Obviously, $R$ crosses $\dot{Q}$ or $\ddot Q$.
If $R$ crossed both $\dot{Q}$ and $\ddot Q$,
then $h(R)\in \{L_{w_{k+1}w_{k+2}\ldots w_{n-1}}(\dot{S}),\linebreak
L_{w_{k+1}w_{k+2}\ldots w_{n-1}}(\ddot{S})\}$ would
separate $S$ and $Q$, which is impossible by 
Corollary \ref{never}. Other\-wise, by the Order-invariance Principle,
$h(R)$ would cross $Q$, in contradiction to 
$h(R)\in {\cal B}(Q)$.
This finishes the induction.\vspace{1,5mm} 

$\dot{Q}=h^{n-1}(Q)$ separates $\dot{S}$ and $\ddot{S}$. 
So the Order-invariance Principle 
shows that $h^{n-2}(Q)$ separates 
$L_{w_{n-1}}(\dot{S})$ and $L_{w_{n-1}}(\ddot{S})$,
then that $h^{n-3}(Q)$ separates the leaves
$L_{w_{n-2}w_{n-1}}(\dot{S})$ and $L_{w_{n-2}w_{n-1}}(\ddot{S})$
etc. Finally, one obtains that
the two leaves $L_{\bf w}(\dot{S})$ and
$L_{\bf w}(\ddot{S})$ are separated by $Q$.
\qquad $\Box$
\section{Some further remarks on the Translation
Principle} To explain where the interest
for the Translation Principle comes from, let us mention
a problem, which originates in the combinatorial description 
of the Mandelbrot set. 
The central concept in the approach by Lau and Schleicher
is the {\em internal address} of a hyperbolic
component (see \cite{laus2,schlei}), which can be given recursively:
\begin{enumerate}
\item[1.] The main component has internal address 1.
\item[2.] If $W$ is a hyperbolic component with internal address
$1\rightarrow n_2\rightarrow n_3\ldots
\rightarrow n_k=m$, then a hyperbolic component visible from $W$ 
and of period $n>m$ has internal address
$1\rightarrow n_2\rightarrow n_3\ldots
\rightarrow n_k=n_{k+1}=n$. 
\end{enumerate}

One easily sees that, by this definition, to each 
hyperbolic component there is assigned
an internal address, but Lau and Schleicher have shown more:
By the procedure defining the internal address, in each step
one fixes a hyperbolic component. 
If one further fixes the internal
angle $\frac{p}{q}$, for which the succeeding hyperbolic component 
is contained in the $\frac{p}{q}$-sublimb of the given one,
then one obtains the {\em angled internal address} (compare 
\cite{laus2}, Definition 6.1). 
Angled internal addresses are
complete, i.e.\ different hyperbolic components have different 
angled internal addresses (see \cite{laus2}, Theorem 9.2).\vspace{1,5mm}

However, not each `abstract internal address' $1\rightarrow n_2\rightarrow
\ldots\rightarrow n_k$ occurs as the internal 
address of a hyperbolic component
(compare \cite{laus2,schlei}), and so
the following problem remains:\vspace{3mm}\\
{\bf Problem (geometric version)} 
{\it Which (angled) internal addresses are 
admissible by a hyperbolic component of the Mandelbrot set?}\vspace{3mm}

In \cite{bk2}, Bandt and the author have given a description
of the (abstract) Mandelbrot set by use of kneading sequences.
This way of description, based on Thurston's ideas in \cite{th},
allows a good insight into the relation between
dynamical properties of Julia sets and properties of the Mandelbrot set
(compare \cite{ke2,ke3}, especially \cite{ke4}).\vspace{1,5mm}

Internal addresses can be turned into kneading sequences 
and vice versa (see \cite{laus2}): Let $1\rightarrow n_2\rightarrow n_3\ldots 
\rightarrow n_k=m$ be the internal address of a hyperbolic component
of period equal to $m$
and let $\alpha\in T$ be given such that the corresponding
external angle lands at
its root. Then the kneading sequence
of $\alpha$ starts with 0, its $m$-th
symbol is $\ast$, and if its initial subword of length 
$n_i$ is ${\bf w}$,
then the initial subword of the length $n_{i+1}$ taken from
the kneading sequence and $\overline{\bf w}$ coincide with exception
of exactly the last symbol.
(This procedure becomes also clear from Lemma \ref{futpas}, and
to find the inverse procedure from the kneading sequence
to the internal address is easy.)\vspace{1,5mm}

In so far, internal addresses are a tool equivalent to kneading sequences,
but more compact and containing the geometric information in a direct
form, and it is not hard to show that the following problem 
can be reduced to find a description of all admissible
internal addresses (compare \cite{ke4}).\vspace{3mm}\\
{\bf Problem (symbolic version)}
{\it Which sequences are admissible by the kneading sequence 
of a point $\alpha\in T$?}\vspace{3mm}

The formulation of this problem, which was
noted by C.\ Penrose (see \cite{pen}) and Bandt and the author
(see \cite{bk}), is on a rather elementary level, and at a
first view the relation to complex quadratic iteration
is surprising. 
But, perhaps, a `naive' research of the structure of kneading sequences
would lead to a structure like ${\cal B}_\ast$ directly.\vspace{1,5mm} 

The idea behind internal addresses is not new.
In one-dimensional real dynamics it appeared
as the concept of {\em cutting times}.
Different characterizations of admissible kneading sequences
(or of the admissible internal addresses) for unimodal maps,
and so for real quadratic maps,
have been given, for example by Collet and Eckmann
\cite{ce}, Milnor and Thurston \cite{mt},
Hofbauer and G.\ Keller \cite{hof,hoko}, and Bruin \cite{brui}.

C.\ Penrose has considered an object more general than the abstract Mandelbrot
set and consisting of 0-1-sequences (see {\cite{pen,pen2}). 
To investigate bifurcations there,
and also to describe Julia sets by symbol sequences, he used the concept
of a `principal nonperiodicity function', a generalized form  
of an internal address. 
(We refer to his viewpoint concerning `non-admissible' kneading sequences
and to his rather deep Theorem 4.2, although the latter is only loosely related to
the present subject.)
Also the `combinatorics of initial subwords' in \cite{bk} touches  
the idea of an internal address (compare \cite{ke4}).\vspace{1,5mm}

One can try to find all admissible internal addresses step by
step, as the above defi\-nition suggests. However, in each step
one has to know which number can be appended
to a given admissible internal address to get a new one. 
In fact, this amounts to finding all hyperbolic components
visible from a given one, but of a greater period.\vspace{1,5mm}

In this, one special case is simpler than the general one:
in \cite{laus2},
Lau and Schleicher have called a hyperbolic component of period $m$   
{\em narrow} 
if there exists no hyperbolic component visible from the given one  
of period less than $m$. Obviously, a hyperbolic component $W$ 
of that type is 
characterized by the fact that the leaf in ${\cal B}_\ast$
corresponding to its root has length
$\frac{1}{2^m-1}$. 

By use
of this, Lau and Schleicher have shown that in the narrow case there exist
hyperbolic components visible from $W$ for all $n>m$ (see Theorem 10.2
in \cite{laus2}),
and also for narrow hyperbolic components, they have proved
the Translation Principle.
(For partial proofs of the Translation and Correspondence Principle
in the non-narrow case, compare Proposition 8.4 and 
Corollary 8.5 in \cite{laus2}.)\vspace{1,5mm}

If $W$ of period $m$ 
fails to be narrow, the situation becomes more complicated,
but the general validity of the Translation Principle shows  
that the existence of a visible hyperbolic component of a given period 
greater than $m$ can be decided by looking at the $\frac{1}{2}$-sublimb
of $W$: 

If this sublimb contains a visible component of period $r$ 
or $m+r$ with
$r=1,2,\ldots,m-1$, then there exist visible components of period
$jm+r$ for all $j=1,2,3,\ldots$. In the other direction,
if there exists a visible component of period
$jm+r$, then let $j_0$ be the least $j$ with this property.
The corresponding visible component lies in the $\frac{1}{2}$-sublimb,
and $j_0$ must be equal to 0 or to 1. So we have
the following
\begin{cor}\label{tracor}
{\rm (Translation Principle for non-angled Internal Addresses)}\\
If the
internal address $1\rightarrow n_2\rightarrow n_3\ldots
\rightarrow n_k=m$ is admissible, then for each $j=2,3,4,\ldots$
the internal address 
$1\rightarrow n_2\rightarrow n_3\ldots
\rightarrow n_k=m\rightarrow jm$
is admissible.
Moreover, for each $r=1,2,\ldots,m-1$, the following holds:
If one of the internal addresses 
$1\rightarrow n_2\rightarrow n_3\ldots
\rightarrow n_k=m\rightarrow jm+r$ with $j=1,2,3,\ldots$ is
admissible, then all these addresses are admissible. \qquad $\Box$
\end{cor}

Clearly, the first part of the theorem is trivial,
and the theorem does not touch the embedding in the plane.
However, it is no problem to say more about this embedding
by use of angled internal addresses. This is left to the reader.\vspace{3mm}

{\bf Acknowledgment} The author would like to thank Dierk Schleicher
for detailed and constructive proposals to improve former versions
of the paper.


\newpage
\appendix
\centerline{\huge\sectfont Errata}

\medskip
\noindent The `Translation' and `Correspondence' Principles we
gave earlier turned out to 
be false in the general case. The aim of this errata is to discuss
which parts of the two statements are incorrect and which parts
remain true.
%
\section{Weaker statements}
\paragraph{Some notations.} We make use of {\bf notations,
statements} and {\bf refe\-ren\-ces} as given in the {\it Stony
Brook IMS Preprint 1997/14}. Beyond this some further notations
will be convenient.\vspace{1,5mm}

By a {\em visibility tree} of a given hyperbolic component $W$ we
understand the tree of hyperbolic components visible from $W$ in a
sublimb of $W$. The visibility tree contained in a
$\frac{p}{q}$-sublimb of $W$ is denoted by ${\mathfrak
Vis}_{\frac{p}{q}}(W)$. Two visibility trees ${\mathfrak
Vis}_{\frac{p_1}{q_1}}(W)$ and ${\mathfrak
Vis}_{\frac{p_2}{q_2}}(W)$ are called {\em equivalent} if they
coincide, including the embedding into the plane, when all periods
in ${\mathfrak Vis}_{\frac{p_1}{q_1}}(W)$ are increased by
$(q_2-q_1)m$.

Further, we denote the hyperbolic component associated with a leaf
$B\in {\cal B}_\ast$ by $W_B$. (Recall that the two end points of
$B$ are the unique external angles of the root of that hyperbolic
component.)
\paragraph{The false statements and their modifications.}
The Translation Principle stated at page 3 does not hold in
general. A counter-example is given by the hyperbolic component
$W=W_{\frac{13}{31}\frac{18}{31}}$ of period $m=5$ (see figure):
${\mathfrak Vis}_{\frac{1}{2}}(W)$ does not contain a hyperbolic
component of period $6(=11-5)$, and beyond this the embedding of
the hyperbolic components of periods 8 and 11 is different in
${\mathfrak Vis}_{\frac{1}{3}}(W)$ and ${\mathfrak
Vis}_{\frac{2}{3}}(W)$.\\ \unitlength1cm

\begin{figure}[h]
\begin{center}
\begin{picture}(12,7)
\put(0.4,6.1){${\mathfrak Vis}_{\frac{3}{4}}(W)$}
\put(2.8,6.1){${\mathfrak Vis}_{\frac{2}{3}}(W)$}
\put(5.2,6.1){${\mathfrak Vis}_{\frac{1}{2}}(W)$}
\put(7.6,6.1){${\mathfrak Vis}_{\frac{1}{3}}(W)$}
\put(10,6.1){${\mathfrak Vis}_{\frac{1}{4}}(W)$}
\put(5.75,-0.4){$W$} \epsfxsize 120mm \epsffile{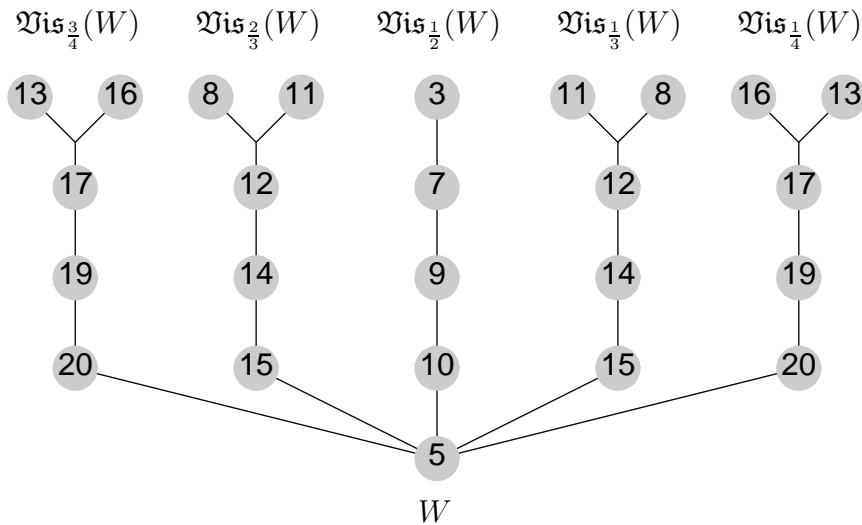}
\end{picture}
\end{center}
\caption{Counter example}
\end{figure}

\noindent The following theorem answers the question of which part
of the Translation Principle can be saved.\vspace{3mm}\\ {\bf
Theorem I.} {\rm (Partial Translation Principle)}\\ {\it Each
visibility tree other than ${\mathfrak Vis}_{\frac{1}{2}}(W)$ for
a hyperbolic component $W$ is equivalent to ${\mathfrak
Vis}_{\frac{1}{3}}(W)$ or to ${\mathfrak
Vis}_{\frac{2}{3}}(W)$.}\vspace{3mm}\\ This errata is mainly
devoted to the proof of Theorem I. Since the Translation Principle
is not valid in general, Corollary 3 must be changed in the
following way:\vspace{3mm}\\ {\bf Corollary I.} {\rm (Partial
Translation Principle for Internal Addresses)}\\ {\it If the
internal address $1\rightarrow n_2\rightarrow n_3\ldots
\rightarrow n_k=m$ is admissible, then for each $j=2,3,4,\ldots$
the internal address $1\rightarrow n_2\rightarrow n_3\ldots
\rightarrow n_k=m\rightarrow jm$ is admissible. Moreover, for each
$r=1,2,\ldots,m-1$, the following holds: If one of the internal
addresses $1\rightarrow n_2\rightarrow n_3\ldots \rightarrow
n_k=m\rightarrow jm+r$ with \fbox{$j=2,3,4,\ldots$} is admissible,
then all these addresses are admissible.\vspace{3mm}}\\ In order
to see that Corollary 3 is false, look at the above
counter-example. The kneading sequence of
$S=\frac{13}{31}\frac{18}{31}$ is $\overline{0100\ast}$ and so
$W=W_S$ has internal address $1\rightarrow 2\rightarrow
4\rightarrow 5$. Therefore, the hyperbolic components of period 11
visible from $W$ have internal address $1\rightarrow 2\rightarrow
4\rightarrow 5\rightarrow 11$.

On the other hand, a hyperbolic component of internal address
$1\rightarrow 2\rightarrow 4\rightarrow 5\rightarrow 6$ does not
exist, providing a counter-example to the statement of Corollary
3. The failure of that internal address can easily be obtained by
checking the kneading sequences of all points in $(T,h)$ of period
6.\vspace{3mm}

The source of all errors in the paper is that Corollary 2 is
false. However, the statements in `3. Proofs' before Corollary 2.
guarantee the following: To visible Fatou components in the
dynamic plane there correspond visible hyperbolic components in
the parameter space, which do not bifurcate directly from $W$.

The problem, hidden by the false Corollary 2, is the existence of
visible hyperbolic components which do not bifurcate directly from
$W$ and are not covered by visible Fatou components. Therefore, in
the general case only the following can be shown:\vspace{3mm}\\
{\bf Proposition I.} {\rm (Partial Correspondence Principle)}\\
{\it Let $W$ be a hyperbolic component of {\em period} $m$, let
$c$ be its center and $A$ be the critical value component of
$p_c$. Then, for relatively prime $r,q$ and a given sublimb of $A$
of denominator $2^{q-1}$ the following holds:

There exists an injective map from the tree of all Fatou
components of {\sc Step} less than $qm$ visible from $A$ and
contained in that sublimb, into ${\mathfrak
Vis}_{\frac{r}{q}}(W)$, which preserves the embedding into the
plane and which turns {\sc Step}'s into periods.\vspace{3mm}}\\
Roughly speaking, the Correspondence Principle transfers `dynamic
translation' to the parameter space; if this transfer is not
complete, the Translation Principle can be incomplete.

The proof of Theorem I will base on a modification of the
Correspondence Principle transferring more `dynamic translation'
to the parameter space than Proposition I (see Proposition II).
\section{Proof of the Partial Translation Principle.}
Let $S\in {\cal B}_\ast$ of period $m$ be fixed, and let ${\bf
v}={\bf v}^S=v_1v_2\ldots v_{m-1}$ and $e=e^S$. Among the boundary
leaves of the critical value gap of ${\cal B}(S)$ there are the
special leaves $R_B$ for visible $B\in {\cal B}_\ast$. We denote
such leaves $R_B$ also by $R_{\frac{p}{q}}$ if the root of $W_B$
has internal angle $\frac{p}{q}$ with respect to $W_S$. ($W_B$
bifurcates directly from $W$ and determines ${\mathfrak
Vis}_{\frac{p}{q}}(W)$.) Note that $R_{\frac{1}{2}}=L_{{\bf
v}e}(S)$ and
$\{h^{m-1}(R_{\frac{1}{3}}),h^{m-1}(R_{\frac{2}{3}})\}=\{L_{e{\bf
v}e}(S),L_{(1-e){\bf v}e}(S)\}$.

By Propositions 1 and 2 and by Lemma 4 one easily gets the
following\vspace{3mm}\\ {\bf Lemma I.} {\it For $p,q$ relatively
prime, the kneading sequence of $R_{\frac{p}{q}}$ is equal to
$({\bf v}e)^{q-1}\overline{{\bf v}(1-e)}$. Moreover, the kneading
sequence of each point $\alpha$ behind $R_{\frac{p}{q}}$ begins
with $({\bf v}e)^{q-2}$.

If $q\geq 3$, then it holds $h^{(q-2)m-1}(R_{\frac{p}{q}})\in
\{L_{e{\bf v}e}(S),L_{(1-e){\bf v}e}(S)\}$, and
$h^{(q-2)m-1}(\alpha)$ lies behind $L_{e{\bf v}e}(S)$ or
$L_{(1-e){\bf v}e}(S)$. \qquad $\blacksquare$\vspace{3mm}}\\ As
already mentioned above there exist visible hyperbolic components
which do not correspond to visible Fatou components in the dynamic
plane. In order to get a better correspondence let us introduce
the concept of semi-visibility. We do this in the lamination
framework.\vspace{3mm}\\ {\bf Definition I.} (semi-visible dynamic
pair)\\ {\it Let $B\in {\cal B}_\ast$ be immediately visible from
$S$. Then a dynamic pair $(S_1,S_2)$ of {\sc Step} $n$ behind
$R_B$ is said to be {\em semi-visible} from $R_B$ if $S_1,S_2\in
{\cal B}$, and if the sequence $\overline{{\bf v}e}$ and the
kneading sequence of $S_1$ have the same initial subwords of
length $n$.}\vspace{3mm}\\ {\bf Remark:} If $B$ has period
$qm;q\geq 2$ and $n<qm$, then visibility of the dynamic pair
$(S_1,S_2)$ implies its semi-visibility (see Corollary 1(ii)).
Moreover, semi-visibility of $(S_1,S_2)$ provides $n<qm$. This
follows from Lemma 2 and the fact that $S_1$ and $S$ are separated
by $R_B$.\vspace{3mm}\\ The statements of Proposition 2 concerning
lengths of leaves show that no iterate of a dynamic pair
$(S_1,S_2)$ of some {\sc Step} $n$ can lie between its members. In
particular, for $S_1=\beta_1\beta_2,S_2= \beta_3\beta_4$ with
$\beta_1<\beta_3<\beta_2<\beta_4$ no iterate of the points
$\beta_1,\beta_2,\beta_3,\beta_4$ lies between $\frac{\beta_1}{2}$
and $\frac{\beta_2}{2}$ or between $\frac{\beta_3}{2}$ and
$\frac{\beta_4}{2}$. This and
$\{h^{n-1}(S_1),h^{n-1}(S_2)\}=\{\dot{S},\ddot{S}\}$ imply the
following statement: \vspace{3mm}\\{\bf Lemma II.} {\it Let $B\in
{\cal B}_\ast$ of period $qm$ be immediately visible from $S$ and
let $(S_1,S_2)$ be a semi-visible dynamic pair of {\sc Step} $n$
behind $R_B$. Then $n<qm$. Moreover, $\overline{{\bf v}e}$ and the
kneading sequence of $S_2$ have the same initial subwords of
length $n$ with the exception of their $n$-th symbols. \qquad
$\blacksquare$\vspace{3mm}}\\ Let us now relate semi-visibility in
the dynamic plane to visibility in the parameter
space.\vspace{3mm}\\ {\bf Proposition II.} (`dynamic'
semi-visibility and `parameter' visibility)\\ {\it Let $B\in {\cal
B}_\ast$ be immediately visible and of period $qm;q\geq 3$. Then
there exists a map $Q:(S_1,S_2)\mapsto Q(S_1,S_2)$ from the set of
all dynamic pairs semi-visible from $R_B$ onto the set of all
visible leaves in ${\cal B}_\ast$ behind $B$ satisfying the
following properties:
\begin{enumerate}
\item[(i)]
The {\sc Step} of $(S_1,S_2)$ and the period of $Q(S_1,S_2)$
coincide.
\item[(ii)]
$Q(S_1,S_2)$ is the only visible leaf in ${\cal B}_\ast$
separating $S_1$ and $S_2$.
\end{enumerate}
Moreover, if $(S_1,S_2)$ is a dynamic pair semi-visible from
$R_B$, then the $(q-3)m$-th iterate of $(S_1,S_2)$ lies behind
$R_{\frac{1}{3}}$ or $R_{\frac{2}{3}}$, and each higher iterate
lying behind $S$ is behind $R_{\frac{1}{2}}$.} \vspace{2,5mm}\\
Proof: Let $(S_1,S_2)$ be semi-visible from $R_B$ of {\sc Step}
$n=jm+k$ with $j<q$ and $0\leq k<m$. Denote the two open intervals
between $S_1$ and $S_2$ each having one end together with $S_1$
and one together with $S_2$ by $I_1$ and $I_2$.

By Proposition 2, these intervals are shorter than
$\frac{1}{2^n}$. Thus, for each $l\leq n$, both $I_1$ and $I_2$
contain at most one periodic point of period $l$ (see above Lemma
1). By Lemma 1, there is one in $I_1$ and one in $I_2$ for $l=n$
since in this case the kneading sequences of $S_1$ and $S_2$ have
different $n$-th symbols, but no one for $l<n$ since the $l$-th
symbols coincide now.\vspace{1,5mm}

Let $Q=Q(S_1,S_2)$ be the leaf connecting the two periodic points
of period $n$ we have obtained. $S_1,S_2\in {\cal B}$ forces that
$Q\in {\cal B}_\ast$, and obviously the kneading sequence of $Q$
is equal to $(\overline{{\bf v}e)^jv_1v_2\ldots v_{k-1}\ast}$.
Moreover, $Q$ is visible by Lemma 2.

We now verify that there is no other leaf in ${\cal B}_\ast$
visible from $S$ and separating $S_1$ and $S_2$. Assuming the
opposite, the period of such a leaf must be greater than $n$ and
so by visibility it must separate $S_1$ and $Q$.\vspace{1,5mm}

Let $r$ be the minimum of periods of all leaves in ${\cal B}_\ast$
separating $S_1$ and $Q$. Then by Lavaurs' Lemma there is a unique
one of period $r$, and $r>n$. Let ${\bf s}$ and ${\bf q}$ be the
kneading sequence of $S_1$ and the kneading sequence just before
$Q$, respectively. Then by Lemma 2, the initial subwords of length
$r-1$ of ${\bf s},{\bf q}$ and $\overline{{\bf v}e}$ coincide, but
the $r$-th symbols of ${\bf s}$ and ${\bf q}$ are different.

Thus from ${\bf s}=({\bf v}e)^j v_1v_2\ldots v_k \overline{{\bf
v}(1-e)}$ and ${\bf q}= (\overline{{\bf v}e)^jv_1v_2\ldots
v_{k-1}v_k}$ one gets $r\geq(j+1)m+k$, implying $({\bf
v}e)^jv_1v_2\ldots v_k{\bf v}= ({\bf v}e)^{j+1}v_1v_2\ldots v_k$.
So, 
we have ${\bf v}=v_{k+2}\ldots v_{m-1}ev_1\ldots v_k$,
which contradicts Lemma 3.\vspace{1,5mm}

It remains to show the last statement of the proposition and that
$Q$ is surjective. For this, let $Q\in {\cal B}_\ast$ behind $B$
be visible. Further, let $n$ be the period of $Q$ and let $l$ be
minimal with the property that $h^l(Q)$ separates $\dot{S}$ and
$\ddot{S}$. Clearly, $l+1\leq n<qm$, and according to Proposition
4(ii), one has $S,\dot{S},\ddot{S}\in {\cal B}(Q)$. Let
$R=h^l(Q)$.

Denote the i-th symbol of $\overline{{\bf v}e}$ by $w_i$ and fix
symbols $s_1,s_2,\ldots, s_l\in \{0,1\}$ such that
$Q=L_{s_1s_2\ldots s_l}(R)$. By induction on $k$ one shows that
for all $k=1,2,\ldots ,l$ the leaves $L_{s_{l-k+1}s_{l-k+2}\ldots
s_l}(\dot{S})$ and $L_{s_{l-k+1}s_{l-k+2}\ldots s_l}(\ddot{S})$
are se\-pa\-rated by $L_{s_{l-k+1}s_{l-k+2}\ldots s_l}(R)=
h^{l-k}(Q)$. (One uses that $L_{s_{l-k+1}s_{l-k+2}\ldots
s_l}(\dot{S}), L_{s_{l-k+1}s_{l-k+2}\ldots
s_l}(\ddot{S})$ and $h^{l-k}(Q)$ do not separate $\dot{S}$ and
$\ddot{S}$.)\vspace{1,5mm}

Let $S_1$ be the longer of the leaves $L_{s_1s_2\ldots
s_l}(\dot{S})$ and $L_{s_1s_2\ldots s_l}(\ddot{S})$, and let $S_2$
be the other one. Clearly, by Propositions 1 and 2 the $(q-3)m$-th
iterates of $S_1,S_2,Q$ lie behind $R_{\frac{1}{3}}$ or
$R_{\frac{2}{3}}$ and their $(q-2)m$-th iterates behind
$R_{\frac{1}{2}}$.

If for some $j>0$ the leaves $h^{(q-2)m+j}(S_1),h^{(q-2)m+j}(S_2)$
were behind $R_{\frac{1}{3}}$ or $R_{\frac{2}{3}}$, then also the
leaf $h^{(q-2)m+j}(Q)$. So by Proposition 1 the word
$w_{(q-2)m+j+1}w_{(q-2)m+j+2}\ldots w_{(q-1)m+j-1}$ would coincide
with ${\bf v}$ and Lemma 3 would imply $j=m$. So $(h^{qm}(S_1),
h^{qm}(S_2))$ would be behind $R_{\frac{1}{2}}$, what is obviously
false.

Therefore, $h^{(q-2)m+j}(S_1),h^{(q-2)m+j}(S_2)$ cannot be behind
$R_{\frac{1}{3}}$ or $R_{\frac{2}{3}}$. In particular, we get that
no iterate of $S_1,S_2$ separates $S_1$ or $S_2$ from $S$.
According to the considerations above Proposition 4, this implies
$S_1,S_2\in {\cal B}$.\vspace{1,5mm}

To show that $(S_1,S_2)$ forms a dynamic pair, assume the
contrary. Then $s_{l-m+2}s_{l-m+3}\ldots s_l={\bf v}$ and
$\{h^{l-m+1}(S_1),h^{l-m+1}(S_2)\}=\{S,R_{\frac{1}{2}}\}$.
Therefore, $h^{l-m+1}(Q)$ separates the boundary chords $S$ and
$R_{\frac{1}{2}}$ of the critical value gap, implying
$w_{l-m+2}w_{l-m+3}\ldots w_l={\bf v}$. Clearly, if $j\leq
(q-2)m$, then $h^j(Q)$ is shorter than $R_{\frac{1}{2}}$. So we
have $l-m+1>(q-2)m$. Hence Lemma 3 provides $l-m+2=(q-1)m+1$,
leading to the contradiction $l+1=qm$.\vspace{1,5mm}

Since $Q$ is visible, $S$ and $S_1$ cannot be separated by a leaf
in ${\cal B}_\ast$ of period less than or equal to $n$. Thus by
Lemma 2 $(S_1,S_2)$ is semi-visible, and the first part of the
proof shows $n=l+1$ and $Q=Q(S_1,S_2)$.

So we have surjectivity of $Q$ and therefore the above also
provides the last statement of the proposition. \qquad
$\blacksquare$\vspace{3mm}\\With Propositions 1 and II and Lemma
III following below, the proof of Theorem I is complete. Namely,
Lemma III relates semi-visibility from different leaves $R_B$, and
Theorem I can be specified as follows:\vspace{3mm}\\ {\bf Theorem
II.} {\it Let $W$ be a hyperbolic component, and let $p_1,q_1$ and
$p_2,q_2$ be relatively prime with $q_1,q_2\geq 3$. If
$h^{(q_1-3)m}(R_{\frac{p_1}{q_1}})=h^{(q_2-3)m}(R_{\frac{p_2}{q_2}})$, 
then the
visibility trees ${\mathfrak Vis}_{\frac{p_1}{q_1}}(W)$ and
${\mathfrak Vis}_{\frac{p_2}{q_2}}(W)$ are equivalent.
\qquad $\blacksquare$\vspace{3mm}}\\ {\bf Lemma III.} {\it Let
$p_1,q_1$ and $p_2,q_2$ be relatively prime with $q_1,q_2\geq 3$
and with
$h^{(q_1-3)m}(R_{\frac{p_1}{q_1}})=h^{(q_2-3)m}(R_{\frac{p_2}{q_2}})$.
Further, let $(S_1,S_2)$ and $(S_3,S_4)$ be dynamic pairs of  {\sc
Step}'s $k<q_1m$ and $l<q_2m$ behind $R_{\frac{p_1}{q_1}}$ and
$R_{\frac{p_2}{q_2}}$, respectively.

If $h^{(q_1-3)m}(S_1)=h^{(q_2-3)m}(S_3)$ and
$h^{(q_1-3)m}(S_2)=h^{(q_2-3)m}(S_4)$, then it holds
$l-k=(q_2-q_1)m$, and $(S_1,S_2)$ is semi-visible iff $(S_3,S_4)$
is.\vspace{2,5mm}}\\ Proof: $l-k=(q_2-q_1)m$ is obvious. So let us
assume that $(S_1,S_2)$ is semi-visible. Then by the last
statement of Proposition II only the $(q_1-3)m$-th iterate of
$(S_1,S_2)$ lies behind $R_{\frac{1}{3}}$ or $R_{\frac{2}{3}}$ and
so only the $(q_2-3)m$-th one of $(S_3,S_4)$ (see also Proposition
1).\vspace{1,5mm}

The first consequence of this fact is that $S_3,S_4\in {\cal
B}_\ast$. Namely, no iterate of $S_3,S_4$ separates $S$ from $S_3$
or $S_4$, and one can argue as in the second part of the proof of
Proposition II. Moreover, together with Proposition 1 one obtains
the following consequence: The leaf
$h^{(q_1-2)m-1}(S_1)=h^{(q_2-2)m-1}(S_3)$ lies behind $L_{e{\bf
v}e}(S)$ or $L_{(1-e){\bf v}e}(S)$, and if some of its iterates is
between $\dot{S}$ and $\ddot{S}$, then again behind one of the
leaves $L_{e{\bf v}e}(S)$ or $L_{(1-e){\bf v}e}(S)$.

Since by Lemma I the kneading sequences of $S_3$ begins with
$({\bf v}e)^{q_2-2}$, the initial subwords of length $l$ of
$\overline{{\bf v}e}$ and the kneading sequence of $S_3$ coincide,
and we are done. \qquad $\blacksquare$\vspace{3mm}\\ Let us finish
by noting the problem in the proof of the false Corollary 2: The
dynamic pair $(S_1,S_2)$ need not be visible because it can lie
between $S_3$ and $S_4$ for some visible dynamic pair $(S_3,S_4)$,
although $S_3$ and $S_4$ are not separated by $S_1,S_2$ (and
$Q$).\newpage

\end{document}